\documentclass[journal]{IEEEtran}

\usepackage{booktabs}
\usepackage{epsfig}
\usepackage{graphicx}
\usepackage{amssymb}
\usepackage{amsfonts}
\usepackage{colortbl}
\usepackage{multirow}
\usepackage{supertabular}
\usepackage{bm}
\usepackage{multicol}
\usepackage{algorithmic}
\usepackage{amsmath}
\usepackage{amsthm}
\usepackage[linesnumbered,ruled,vlined]{algorithm2e}
\usepackage{mathtools}
\usepackage[dvipsnames,table]{xcolor}
\usepackage{subfigure}
\usepackage{paralist}
\usepackage{color}
\usepackage{tabu}
\usepackage{rotating}
\usepackage{pifont}
\usepackage{cite}
\usepackage{threeparttable}
\usepackage{moresize}
\usepackage{cancel}
\usepackage{pgf-pie}
\usepackage{adjustbox}
\usepackage[margin=2cm]{geometry}
\usepackage{authblk}
\usepackage{url}
\usepackage[numbers,sort&compress]{natbib}
\usepackage[graphicx]{realboxes}
\usepackage{subcaption}
\usepackage[colorlinks=true,linkcolor=black,anchorcolor=black,citecolor=black,filecolor=black,menucolor=black,runcolor=black,urlcolor=blue]{hyperref}

\newcommand{\CM}{\ding{51}}
\newcommand{\XM}{\ding{53}}
\newcolumntype{L}{>{\raggedleft}p{0.14\textwidth}}
\newcolumntype{R}{p{0.8\textwidth}}

\newtheorem*{example*}{Example}

\renewcommand{\vec}{\mathbf}

\title{\textbf{GNBG-Generated Test Suite for Box-Constrained Numerical Global Optimization}}
\author{Amir H. Gandomi$^{\ast\dagger}$, Danial~Yazdani$^\ast$, Mohammad~Nabi~Omidvar, and Kalyanmoy~Deb
    
\thanks{Amir H. Gandomi is with the Faculty of Engineering \& Information Technology, University of Technology Sydney, Ultimo 2007, Australia. 
He is also with the University Research and Innovation Center (EKIK), Obuda University, Budapest 1034, Hungary.
(e-mail: Gandomi@uts.edu.au)}

\thanks{Danial Yazdani is with the Faculty of Engineering \& Information Technology, University of Technology Sydney, Ultimo 2007, Australia. (e-mail: danial.yazdani@gmail.com)
}

\thanks{Mohammad Nabi Omidvar is with the School of Computing, University of Leeds, and Leeds University Business School, Leeds LS2 9JT, United Kingdom. (e-mail: m.n.omidvar@leeds.ac.uk)}

\thanks{Kalyanmoy Deb is with the BEACON Center, Michigan State University, East Lansing, MI, 48824, USA. (e-mail: kdeb@egr.msu.edu)}

\thanks{\emph{$^\ast$Amir H. Gandomi and Danial Yazdani contributed equally to this work.}}
\thanks{\emph{$^\dagger$Corresponding author: Amir H. Gandomi}}

}

\begin{document}

\maketitle

\begin{abstract}

This document introduces a set of 24 box-constrained numerical global optimization problem instances, systematically constructed using the Generalized Numerical Benchmark Generator (GNBG). 
These instances cover a broad spectrum of problem features, including varying degrees of modality, ruggedness, symmetry, conditioning, variable interaction structures, basin linearity, and deceptiveness. 
Purposefully designed, this test suite offers varying difficulty levels and problem characteristics, facilitating rigorous evaluation and comparative analysis of optimization algorithms. 
By presenting these problems, we aim to provide researchers with a structured platform to assess the strengths and weaknesses of their algorithms against challenges with known, controlled characteristics. 
For reproducibility, the MATLAB source code for this test suite is publicly available and can be accessed at~\cite{yazdani2023GNBGinstances}.

\end{abstract}

 
 \section{Introduction}

Optimization algorithms have been the subject of intense research and development over the past decades, finding applications in a variety of domains. One fundamental research question revolves around determining an algorithm's effectiveness on problems with specific characteristics. While theoretical analyses provide insights, they often become challenging for complex algorithms and intricate problem instances. Hence, empirical evaluation, typically executed by solving a predefined set of benchmark problem instances, becomes essential~\cite{beiranvand2017best}.

This document presents a collection of 24 problem instances tailored for benchmarking algorithms in box-constrained numerical global optimization. 
Within this context, benchmarking refers to the comparison of best solutions found by various algorithms, assessed through specific performance indicators~\cite{hansen2021coco}. The global optimization problems are widespread in fields like mathematics and engineering~\cite{mei2021structural}. Employing appropriate benchmark test suites in this domain is not just an academic exercise; it lays the foundation for advancements in more complex optimization problems, including dynamic~\cite{yazdani2021DOPsurveyPartA,yazdani2021DOPsurveyPartB}, constrained~\cite{kusakci2012constrained}, large-scale~\cite{omidvar2021reviewA,omidvar2021reviewB}, niching~\cite{li2016seeking}, and multi-objective optimization~\cite{saini2021multi}.

The test suite's problem instances are derived from the Generalized Numerical Benchmark Generator (GNBG)\footnote{The MATLAB source code for the GNBG problem instance generator is available at~\cite{yazdani2023GNBGgenerator}. Users can employ this code to generate custom problem instances as per their requirements.}, a tool crafted for the systematic evaluation of global optimization algorithms. 
GNBG utilizes a singular parametric baseline function, and by adjusting its parameter settings, it can produce a wide array of problem instances with controllable characteristics and varying degrees of difficulty. 
By manipulating various parameters within GNBG, researchers gain the ability to tailor the properties (individual or any desired combinations) of the generated problem instances, including:
\begin{itemize}
\item Modality: GNBG offers the versatility of generating a diverse range of problem instances, from those that feature smooth, unimodal search spaces to instances characterized by highly multimodal and rugged landscapes. 
This adaptability allows researchers to assess how well optimization algorithms navigate different types of terrain, thereby providing a more comprehensive evaluation of their capabilities.

\item Local Optima Characteristics: GNBG constructs its search space through the integration of multiple independent components, each having its own `basin of attraction'--essentially, a zone where solutions tend to converge. 
Users can configure various aspects of these components, such as their locations, optimum values, and morphological features. 
This high level of control extends to the characteristics of any local optima within these basins, allowing for customization of their number, size, width, depth, and shape.

\item Gradient Characteristics: GNBG allows users to control not just the steepness of the components but also the specific rate of change or curvature of their basins. 
Users have the flexibility to define these characteristics on a per-component basis, with options ranging from highly sub-linear to super-linear rates of change.

\item Variable Interaction Structures: GNBG offers a nuanced control over the interactions between variables within generated problem instances. 
Users can then apply customization rotation matrices to configure variable interactions, allowing for a wide range of interaction structures from fully separable to fully-connected non-separable. 
GNBG allows users to set the strength of these interactions, ensuring a more robust evaluation of algorithms. 
Additionally, different regions of the search space can have distinct variable interaction structures, as each component possesses its own localized variable interaction pattern.

\item Conditioning: GNBG provides users with the ability to generate components with a wide range of condition numbers, spanning from well-conditioned to severely ill-conditioned components. 
By independently stretching each component along each dimension, users have control over the condition number, allowing them to simulate the challenges posed by ill-conditioning.

\item Symmetry: GNBG affords the flexibility to generate both symmetric and highly asymmetric problem instances. 
This is achieved by allowing the strategic distribution of components with varied morphological characteristics across the search space, thereby influencing the symmetry of the resulting problem instances. 
Furthermore, GNBG provides the capability to generate components with asymmetric basins of attraction. 
The ability to introduce asymmetry serves an important role in unbiased evaluations; symmetric search spaces can potentially favor algorithms that rely on Gaussian distribution-based operators~\cite{hansen2009real}. 
Therefore, incorporating asymmetric elements allows for more impartial assessments of algorithmic performance.

\item Deceptiveness: GNBG affords users the flexibility to introduce various degrees of deception into problem instances. 
By manipulating the size, location, and depth of components, one can engineer situations that present specific challenges for optimization algorithms. 
Scenarios can be constructed where, for instance, the global optimum is a narrow peak obscured within the basin of a wider local optimum, or where there is substantial separation between high-quality local optima and the global optimum.
This enables researchers to scrutinize how well algorithms can navigate misleading or complex landscapes, thereby yielding insights into their robustness and efficacy.

\item Scalability: All problem instances generated by GNBG are scalable with respect to dimensionality. 
\end{itemize}

The search space in GNBG is a composite landscape formed by aggregating multiple distinct components, each characterized by its unique basin of attraction. 
These components, which must all have the same dimensionality, contribute individual challenges and complexities to the overall search space.
The baseline function of GNBG is shown in Equation~\eqref{eq:irGMPB}, where, $f(\cdot)$ represents the GNBG function aimed to be minimized. 
In this equation, $d$ indicates the number of dimensions, and $\mathbb{X}$ denotes the $d$-dimensional search space where $\vec{x}$ is a candidate solution. 
The search space is confined within the bounds $l_i$ and $u_i$ for each dimension $i$. 
The parameter $o$ enumerates the number of components, each contributing to the complexity of the problem with its own set of parameters such as $\sigma_k$, $\vec{m}_k$, $\mathbf{H}_{k}$, $\mathbf{R}_k$, and $\lambda_k$.

\begin{figure*}[th!] 
\begin{align}
\label{eq:irGMPB}
\mathrm{Minimize\;\;:} & f(\vec{x})= \min_{k\in\{1,\dots,o\}}\left\{ \sigma_{k} + \bigg(\mathbb{T}_k\Big((\mathbf{R}_{k}(\vec{x}-\vec{m}_k))^\top\Big)  \mathbf{H}_{k} \mathbb{T}_k\Big(\mathbf{R}_{k}(\vec{x}-\vec{m}_k)\Big)\bigg)^{\lambda_{k}} \right\}\\
\mathrm{Subject~to:} & \;\;\vec{x}\in \mathbb{X}:\;\mathbb{X}=\{\vec{x}\;|\;l_i \leq x_i \leq u_i\}, \; i \in\{1,2,\ldots,d\},\nonumber
\end{align}
\end{figure*}

\begin{figure*}[th!] 
\begin{align}
\label{eq:ir}
    a_j \mapsto
    \begin{dcases}
    \exp{\bigg(\log(a_j)+\mu_{k,1}\Big(\sin{\big(\omega_{k,1}\log(a_j)\big)}+\sin{\big(\omega_{k,2}\log(a_j)\big)} \Big)\bigg)} & \text{if   } a_j>0 \\
    0 & \text{if   }a_j=0\\
    -\exp{\bigg(\log(|a_j|)+\mu_{k,2}\Big(\sin{\big(\omega_{k,3}\log(|a_j|)\big)}+\sin{\big(\omega_{k,4}\log(|a_j|)\big)} \Big)\bigg)} & \text{if   } a_j<0  
\end{dcases}
\end{align}
\end{figure*}

\begin{algorithm*}[!tp]
\footnotesize
$\mathbf{R}_{k}=\mathbf{I}_{d \times d}$\;
\For{$p=1$ to $d-1$}{\label{algline:Loop1}
\For{$q=p+1$ to $d$}{\label{algline:Loop2}
\If{$\bm\Theta_{k}(p,q)\neq 0$}{
$\mathbf{G}=\mathbf{I}_{d \times d}$\;
$\mathbf{G}(p,p)= \cos\left(\bm\Theta_{k}(p,q)\right)$;~\tcp*[h]{$\bm\Theta_{k}(a,b)$ and $\mathbf{G}(a,b)$ are the elements at $a$th row and $b$th column of matrices $\bm\Theta_{k}$ and $\mathbf{G}$, respectively.}\\
$\mathbf{G}(q,q) = \cos\left(\bm\Theta_{k}(p,q)\right)$\;
$\mathbf{G}(p,q) =  -\sin\left(\bm\Theta_{k}(p,q)\right)$\;
$\mathbf{G}(q,p) =  \sin\left(\bm\Theta_{k}(p,q)\right)$\;
$\mathbf{R}_{k} = \mathbf{R}_{k} \times \mathbf{G}$;~\tcp*[h]{$\mathbf{G}$ is the Givens rotation matrix for $x_{p}-x_{q}$ plane based on $\bm\Theta_{k}(p,q)$.}\\
}}}
Return $\mathbf{R}_{k}$\;
\caption{Pseudo code for calculating the rotation matrix $\mathbf{R}_{k}$ based on $\bm\Theta_k$.\\
{\footnotesize\textbf{Input}: $d$ and $\bm\Theta_k$\\
\textbf{Output}: {$\mathbf{R}_{k}$}}
}
\label{alg:RotationControlled}
\end{algorithm*}

For the $k$th component, $\vec{m}_k$ defines its center, $\sigma_k$ specifies its minimum value (i.e., $f(\vec{m}_k)$), and $\mathbf{H}_{k}$ is a $d \times d$ diagonal matrix whose principal diagonal elements affect the heights of the basin associated with the $k$th component across different dimensions. 
Further, $\mathbf{R}_k$ serves as the rotation matrix for the $k$th component, and $\lambda_k$ quantifies the linearity degree in its basin of attraction. The $\min(\cdot)$ function delineates the basins of attraction for each component.
Finally, $\mathbb{T}_{k}(\vec{a})\mapsto \vec{b}$ is a non-linear transformation function~\cite{hansen2009real} that introduces additional complexities to the basin of each component.
$\mathbb{T}_{k}(\vec{a})\mapsto \vec{b}$ is formulated as Equation~\eqref{eq:ir}.
Based on this equation, for the $k$th component, this transformation is guided by parameters $\bm\mu_k$ and $\bm\omega_k$, which define the symmetry and morphology of local optima on the basin of the $k$th component.
Table~\ref{tab:GNBGparameters} summarizes the parameters of GNBG.

Algorithm~\ref{alg:RotationControlled} generates a rotation matrix, denoted as $\mathbf{R}_k$, that rotates the projection of $\vec{x}$ within the basin of attraction of the $k$th component based on a given matrix $\bm\Theta_k$.
By setting $\bm\Theta_{k}(p,q)$ to values that are not multiples of $\frac{\pi}{2}$ (i.e., non-axis-overlapped values), we establish variable interactions between $p$ and $q$.
The strength of this interaction is determined by how far $\bm\Theta_{k}(p,q)$ deviates from multiples of $\frac{\pi}{2}$.
For example, a $\bm\Theta_{k}(p,q)$ value of $\frac{\pi}{4}$ would induce a strong variable interaction, while a value close to $\frac{\pi}{180}$ would result in weaker interactions.
Setting $\bm\Theta_{k}(p,q)$ to zero ensures that $\mathbf{R}_k$ does not modify the interaction between $p$th and $q$th variables.
In GNBG, Algorithm~\ref{alg:RotationControlled} is used to configure $\mathbf{R}_k$ to generate controlled variable interaction structures based on the values included in $\bm\Theta_k$.

\begin{table*}[t!] 
\footnotesize
\centering
  \caption{Summary of the parameters, functions, and notations used in GNBG.}
  \label{tab:GNBGparameters}
  \begin{tabular}{cp{15cm}}
\midrule
Symbol & Description\\
\midrule
$f(\cdot)$ & GNBG's baseline function. \\
\rowcolor{gray!20}$d$ & Represents the number of components in the search space.\\
$\mathbb{X}$ & $d$-dimensional search space.\\
\rowcolor{gray!20}$\vec{x}$ & A $d$-dimensional solution $(x_1,x_2,\ldots,x_d)$ in the search space $\mathbb{X}$. \\
$l_i$ &  Lower bound of search range in the $i$th dimension.\\
\rowcolor{gray!20}$u_i$ &  Upper bound of search range in the $i$th dimension.\\
$o$ & Number of components in the search space.\\
\rowcolor{gray!20}$\mathrm{min}(\cdot)$ & Defines the basin of attraction of each component.\\
$\vec{m}_{k}$ & Minimum position of the $k$th component.\\
\rowcolor{gray!20}$\sigma_{k}$ & Minimum value of the $k$th component, i.e., $f(\vec{m}_{k})=\sigma_{k}$.\\
$\mathbf{H}_{k}$ & It is a $d \times d$ diagonal matrix, i.e., $\mathbf{H}_{k} = \mathrm{diag}(h_1,h_2, \ldots,h_d)\in \mathbb{R}^{d \times d}$, where $h_i = \mathbf{H}_{k}(i,i)$.
The principal diagonal elements of $\mathbf{H}_{k}$ are scaling factors that influence the heights of the basin associated with the $k$th component across different dimensions. 
Moreover, this matrix directly affect the condition number of the $k$th component. 
\\
\rowcolor{gray!20}$\mathbf{R}_{k}$& It is a $d \times d$ orthogonal matrix used for rotating the $k$th component.\\
$\bm\Theta_{k}$ & It is a $d \times d$ matrix, whose elements are utilized to compute the rotation matrix $\mathbf{R}_{k}$ by Algorithm~\ref{alg:RotationControlled}. 
The elements on and below the principal diagonal of $\bm\Theta_{k}$ are zero. 
An element located at the $p$th row and $q$th column of $\bm\Theta_{k}$, denoted as $\bm\Theta_{k}(p,q)$, specifies the rotation angle for the plane $x_{p}–x_{q}$, where $p < q$.
In essence, $\bm\Theta_{k}(p,q)$ governs the interaction between variables $x_p,x_q\in \vec{x}$.\\
\rowcolor{gray!20}$\mathfrak{p}_k$& In cases where $\bm\Theta_{k}$ is randomly generated, $0 \leq \mathfrak{p}_k \leq 1$ controls the random generation of elements above the principal diagonal in the $\bm\Theta_{k}$ matrix.
For each $\bm\Theta_{k}(p,q)$ where $p<q$, a random number is drawn from a uniform distribution.
If this number is less than or equal to $\mathfrak{p}_k$,  $\bm\Theta_{k}(p,q)$ is set to zero.
Otherwise, $\bm\Theta_{k}(p,q)$ is assigned a predefined or a random angle.\\
$\lambda_{k}$ & It is a positive constant that affects the rate at which the basin of the $k$th component increases. 
The specific pattern can range from super-linear ($\lambda_{k}>0.5$) to linear ($\lambda_{k}=0.5$) to sub-linear ($0<\lambda_{k}<0.5$).\\
\rowcolor{gray!20}$\mathbb{T}(\cdot)$ & An element-wise non-linear transformation that plays a role in controlling the modality, irregularity, roughness, and symmetry of each component. 
 The specific characteristics of the $k$th component are determined by the values of $\bm\mu_{k}$ and $\bm\omega_{k}$ within the transformation function.\\
$\bm\mu_{k}$ &   The vector $\bm\mu_{k}$ consists of two elements, denoted as $\bm\mu_k = (\mu_{k,1},\mu_{k,2})$. 
These elements play a crucial role in determining the depth of the local optima within the basin of the $k$th component. 
By assigning different values to $\mu_{k,1}$ and $\mu_{k,2}$, an asymmetry is introduced into the basin of the $k$th component.\\
\rowcolor{gray!20}$\bm\omega_{k}$ &  The vector $\bm\omega_{k}$ consists of four elements, denoted as $\bm\omega_{k} = (\omega_{k,1},\omega_{k,2},\omega_{k,3},\omega_{k,4})$. 
Together with the $\bm\mu_{k}$ values, these elements play a significant role in shaping the characteristics of the "basin optima" within the basin of the $k$th component. 
The values of $\omega_{k,1},\omega_{k,2},\omega_{k,3}$, and $\omega_{k,4}$ contribute to determining the number and width of local optima within the basin of the $k$th component. 
Furthermore, differences between these elements impact the symmetry of the basin of the $k$th component.\\
   \bottomrule
  \end{tabular}
 \end{table*}

In the following section, we present 24 carefully crafted problem instances. 
Utilizing this test suite allows for a systematic exploration of an algorithm's strengths and vulnerabilities across diverse, controlled conditions and challenges. 
Such thorough evaluations are crucial for honing and refining optimization algorithms.

\section{Problem Instances Generated by GNBG}
\label{sec:instance}

In this section, we introduce a set of 24 problem instances, denoted as $f_{1}$ through $f_{24}$, all generated using GNBG. 
These instances are specifically designed to span a wide array of characteristics, establishing a robust test suite for optimization algorithm evaluation.
The problems can be categorized as:

\begin{itemize}
    \item $f_{1}$ to $f_{6}$: Unimodal instances,
    \item $f_{7}$ to $f_{15}$: Multimodal instances with a single component, and
    \item $f_{16}$ to $f_{24}$: Multimodal instances with multiple components.
\end{itemize}
All instances operate in a 30-dimensional solution space and have their search boundaries set within $[-100, 100]^d$.

The primary aim of this test suite is to facilitate in-depth analysis of optimization algorithm behavior and gain a better understanding of their strengths and weaknesses across different problem characteristics and challenges. 
An overview of the 24 problem instances is presented in Table~\ref{tab:SuiteOverview}. 
The MATLAB source code for these problem instances is available for access~\cite{yazdani2023GNBGinstances}. 
The ensuing subsections provide an in-depth overview of each problem instance, elaborating on their defining characteristics and settings.

\begin{table}[t!] 
\footnotesize
\centering
  \caption{An overview of characteristics of the presented test suite containing 24 problem instances generated by GNBG for evaluating optimization algorithms. }  \label{tab:SuiteOverview}
\Rotatebox{90}{
\resizebox{1.22\textwidth}{!}{%
 \begin{threeparttable}
  \begin{tabular}{lcccccccccccccccccccccccccc}
    \toprule
   \multirow{2}{*}{Characteristic} & \multicolumn{24}{c}{Problem instances}\\ 
    \cmidrule{2-25} 
   & $f_{1}$ & $f_{2}$ & $f_{3}$ & $f_{4}$ & $f_{5}$ & $f_{6}$ & $f_{7}$ & $f_{8}$ & $f_{9}$ & $f_{10}$ & $f_{11}$ & $f_{12}$ & $f_{13}$ & $f_{14}$ & $f_{15}$ & $f_{16}$ & $f_{17}$ & $f_{18}$ & $f_{19}$ & $f_{20}$ & $f_{21}$ & $f_{22}$ & $f_{23}$& $f_{24}$  \\
\midrule
  Modality& \multicolumn{6}{c}{Unimodal}&   \multicolumn{9}{c}{\cellcolor{gray!15}Multimodal with single component}& \multicolumn{9}{c}{Multimodal with multiple competing components}\\
\midrule
   Basin local optima\tnote{$\bullet$}& \XM & \XM & \XM & \XM & \XM & \XM & \CM & \CM & \CM & \CM & \CM & \CM & \CM & \CM & \CM & \XM & \XM  & \CM & \CM & \CM & \CM  & \CM & \CM & \CM \\
\midrule
   Separability & $\mathcal{F}$\tnote{$\bm \dagger$} & $\mathcal{F}$ & $\mathcal{F}$ & $\mathcal{N}$\tnote{$\bm \dagger$} & $\mathcal{N}$ & $\mathcal{N}$ & $\mathcal{F}$ & $\mathcal{F}$ & $\mathcal{F}$ & $\mathcal{F}$ & $\mathcal{N}$ & $\mathcal{P}$\tnote{$\bm \dagger$} & $\mathcal{N}$ & $\mathcal{N}$ & $\mathcal{N}$ & $\mathcal{N}$ & $\mathcal{N}$ & $\mathcal{N}$ &  $\mathcal{N}$ & $\mathcal{N}$ & $\mathcal{N}$ & $\mathcal{N}$ & $\mathcal{N}$& $\mathcal{N}$ \\
   \midrule
   Varying variable interactions &  \XM & \XM & \XM & \XM & \XM & \XM & \XM & \XM & \XM & \XM & \XM & \XM & \XM & \XM & \XM & \XM & \CM & \CM & \CM & \CM & \CM & \CM & \CM & \CM  \\
     \midrule
   Symmetry &  $\mathcal{S}$\tnote{$\bm\ddagger$} &  $\mathcal{S}$ &  $\mathcal{S}$ &  $\mathcal{S}$ &  $\mathcal{S}$ &  $\mathcal{S}$ &  $\mathcal{S}$ &  $\mathcal{S}$ &  $\mathcal{S}$ &  $\mathcal{A}$\tnote{$\bm\ddagger$} &  $\mathcal{A}$ &  $\mathcal{A}$ &   $\mathcal{S}$  & $\mathcal{A}$ &  $\mathcal{S}$ &  $\mathcal{A}$ &  $\mathcal{A}$ &   $\mathcal{A}$ &  $\mathcal{A}$ &  $\mathcal{A}$ &  $\mathcal{A}$ &  $\mathcal{A}$ &  $\mathcal{A}$ &  $\mathcal{A}$  \\
     \midrule
     Ill-conditioning &  \XM & \XM & \CM & \XM & \CM & \CM & \XM & \XM & \XM & \XM & \XM & \XM & \XM & \CM & \CM & \XM & \CM & \XM & \XM & \XM & \XM & \XM & \XM & \CM \\
      \midrule
    Basin linearity &  $\mathcal{E}$\tnote{$\bm\ast$} &  $\mathcal{L}$\tnote{$\bm\ast$} &  $\mathcal{E}$ &  $\mathcal{E}$ &  $\mathcal{L}$ &  $\mathcal{L}$ &  $\mathcal{E}$ &  $\mathcal{E}$ &  $\mathcal{E}$ &  $\mathcal{E}$ &  $\mathcal{E}$ &  $\mathcal{E}$ &   $\mathcal{E}$  & $\mathcal{E}$ &  $\mathcal{L}$ &  $\mathcal{E}$ &  $\mathcal{E}$ &  $\mathcal{E}$ &  $\mathcal{E}$ &  $\mathcal{L}$ & $\mathrm{L}$\tnote{$\bm\ast$} & $\mathcal{E}$ &  $\mathcal{L}$   &  $\mathcal{L}$  \\
      \midrule
     Deceptive &  \XM & \XM & \XM & \XM & \XM & \XM & \XM & \XM & \XM & \XM & \XM & \XM & \XM & \XM & \XM & \CM & \CM & \CM & \CM & \CM & \CM & \CM & \XM  & \CM \\
    \bottomrule
  \end{tabular}
 \begin{tablenotes}
\item[$\bullet$] {Existence of local optima within the basin of each component.}
\item[$\bm\dagger$] {$\mathcal{F}$, $\mathcal{N}$ and  $\mathcal{P}$ stand for fully separable, non-separable, and partially separable, respectively.}
\item[$\bm\ddagger$] {$\mathcal{S}$ and $\mathcal{A}$ stand for symmetric and asymmetric, respectively.}
\item[$\bm\ast$] {$\mathcal{E}$, $\mathcal{L}$, and $\mathrm{L}$ stand for super-linear, sub-linear, and linear, respectively.}
\end{tablenotes}
\end{threeparttable}}}
 \end{table}

\subsection{Unimodal problem instances}
\label{sec:sec:unimodalScenarios}
This subsection introduces a set of unimodal problem instances, each exhibiting a single component in the landscape (i.e., $o=1$). 
These functions are defined in a 30-dimensional space ($d=30$) with parameters $\bm\mu=(0,0)$ and $\bm\omega=(0,0,0,0)$. 
The minimum position $\vec{m}$ (global optimum position) is randomly selected from the range $[-80,80]^d$, and $\sigma$ (global optimum value) is randomly selected from the range $[-1,200,0]$.
Note that the random seed used for generating the minimum position $\vec{m}$ remains constant for each problem instance with a single component. This uniform seed ensures consistent and stable performance measurements across all runs.

\subsubsection{$f_1$}

The first instance, labeled as $f_1$, is constructed by setting $\lambda=1$ with $\mathbf{H}=\mathbf{R}=\mathbf{I}_{d \times d}$. 
It mirrors characteristics found in the widely recognized Sphere function. 
This problem is symmetric, regular, smooth, well-conditioned, and fully separable, which facilitates dimension-wise exploitation.
Researchers can employ this instance to examine the convergence speed of optimization algorithms. 
Figure~\ref{fig:f1} visualizes a 2-dimensional representation of $f_1$.

\subsubsection{$f_2$}

The configuration for $f_2$ largely parallels that of $f_1$. 
The distinguishing factor is the setting of $\lambda=0.05$ for this instance. 
Much like $f_1$, $f_2$ is symmetric, regular, well-conditioned, and separable. 
However, in distinction from $f_1$, $f_2$ boasts a sub-linear basin, leading to a pronounced and constricted morphology around its optimal point. 
This attribute tests the optimization algorithms' precision in pinpointing the global optimum within the confined region, rendering it instrumental for gauging algorithms' exploitation efficacy. 
Given that $f_2$ is devoid of other marked challenging traits, it acts as a specific benchmark to evaluate the precision of algorithmic exploitation. 
Figure~\ref{fig:f2} illustrates a 2-dimensional projection of $f_2$.

\subsubsection{$f_3$}

For the configuration of $f_3$, we set $\lambda=1$. 
To initialize the principal diagonal of the matrix $\mathbf{H}$, a linearly spaced set of values ranging from 0.1 to $10^6$ is first generated, inclusive of both endpoints. 
These values are then randomly permuted to assign to the elements of $\mathbf{H}$'s principal diagonal, ensuring a uniform but non-sequential distribution across this range.
Moreover, we use $\mathbf{R}=\mathbf{I}_{d \times d}$ for this problem instance. 

Analogous to $f_1$, this instance retains the traits of smoothness, regularity, and full separability. 
Yet, $f_3$ stands out as a notably ill-conditioned unimodal problem instance with a condition number of $10^7$. 
It can be perceived as an elongated variant of $f_1$, allowing a comparative analysis of algorithmic performance in the face of ill-conditioned challenges. 
The intricacy in $f_3$ arises from the pronounced disparity in the magnitudes of the elements on the principal diagonal of $\mathbf{H}$, leading to the observed ill-conditioning. 
Comparing the results obtained on $f_1$ and $f_3$ provides insights into an algorithm's capability in managing ill-conditioned scenarios and helps discern any performance variances between the two instances.
Figure~\ref{fig:f3} offers a visual representation of $f_3$ in a 2-dimensional space.

\subsubsection{$f_4$}

The first three problem instances ($f_1$ to $f_3$) are fully separable, which facilitates dimension-wise optimization. 
Notably, both $f_1$ and $f_2$ possess condition numbers of unity, rendering them rotation-invariant. 
This invariant property means that the component's orientation remains unchanged even when rotated around its centroid. 
To instill rotation-dependency in these instances, their condition numbers must exceed unity.

For the configuration of $f_4$, we set $\lambda=1$, with every element on the principal diagonal of $\mathbf{H}$ being uniformly randomized within the range [1, 10]. 
Such a setup births a unimodal problem instance characterized by a condition number lying between one and ten. 
Consequently, $f_4$ exhibits rotation-dependency while still being moderately well-conditioned. 
The associated rotation matrix $\mathbf{R}$ gets formulated via Algorithm~\ref{alg:RotationControlled}, applying $\mathfrak{p}=1$, and each angle in $\bm\Theta$ is uniformly drawn from the interval $(-\pi,\pi)$.

Distinctly, $f_4$ emerges as a unimodal, regular, smooth, and non-separable instance. 
While its condition number is larger than one, it is purposely not framed as a severely ill-conditioned problem. 
The objective here is to gauge the influence of the non-separable variable interaction framework on the behavior of optimization algorithms without the overpowering effect of extreme ill-conditioning. 
Drawing a parallel between $f_1$ and $f_4$ will allow for a targeted assessment, isolating the effect of non-separability on algorithmic performance. Figure~\ref{fig:f4} graphically represents a 2-dimensional projection of $f_4$.

%

\subsubsection{$f_5$}

For $f_5$, we set $\lambda=0.05$. 
The elements of the principal diagonal of the matrix $\mathbf{H}$ are determined using a method analogous to that employed in $f_3$: we first generate a linearly spaced set of values ranging from 0.1 to $10^6$, inclusive of both endpoints. 
These values are then randomly permuted to assign to the principal diagonal of $\mathbf{H}$.
In $f_5$, each variable $i$, interacts with the succeeding variable, $(i+1)$, resulting in a chain-like variable interaction structure that defines a minimally connected, non-separable problem instance. 
The rotation matrix $\mathbf{R}$ is computed using Algorithm~\ref{alg:RotationControlled}, with the rotation angles for each plane being randomly determined from the interval $(-\pi,\pi)$.

Characterizing $f_5$ is its non-separability and an imposing condition number of $10^7$. 
Moreover, it exhibits a profoundly sub-linear basin. 
The synergy of ill-conditioning, sub-linear basin, and chain-like variable interaction crafts a distinctive landscape: a sharply defined, rotated valley. 
Such a landscape typically presents an arduous challenge for many optimization algorithms, especially concerning fine-grained exploitation and rapid convergence.

\begin{figure*}[!t]
\centering
\begin{tabular}{ccc}
    \subfigure[{\scriptsize $f_1$.}]{\includegraphics[width=0.30\linewidth]{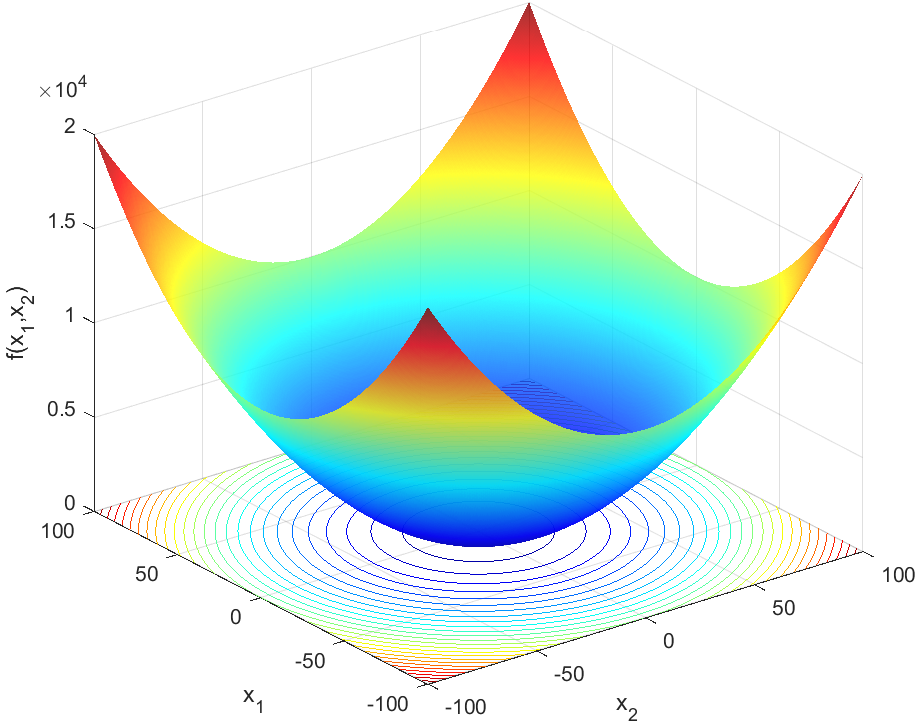}\label{fig:f1}}
&
     \subfigure[{\scriptsize $f_2$.}]{\includegraphics[width=0.30\linewidth]{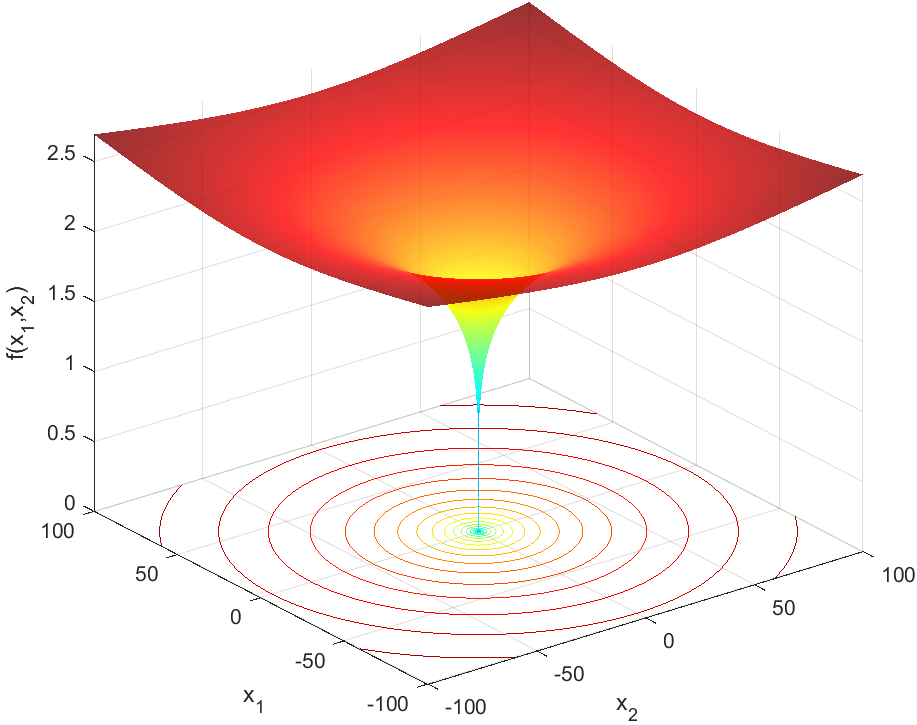}\label{fig:f2}}
&
\subfigure[{\scriptsize $f_3$.}]{\includegraphics[width=0.30\linewidth]{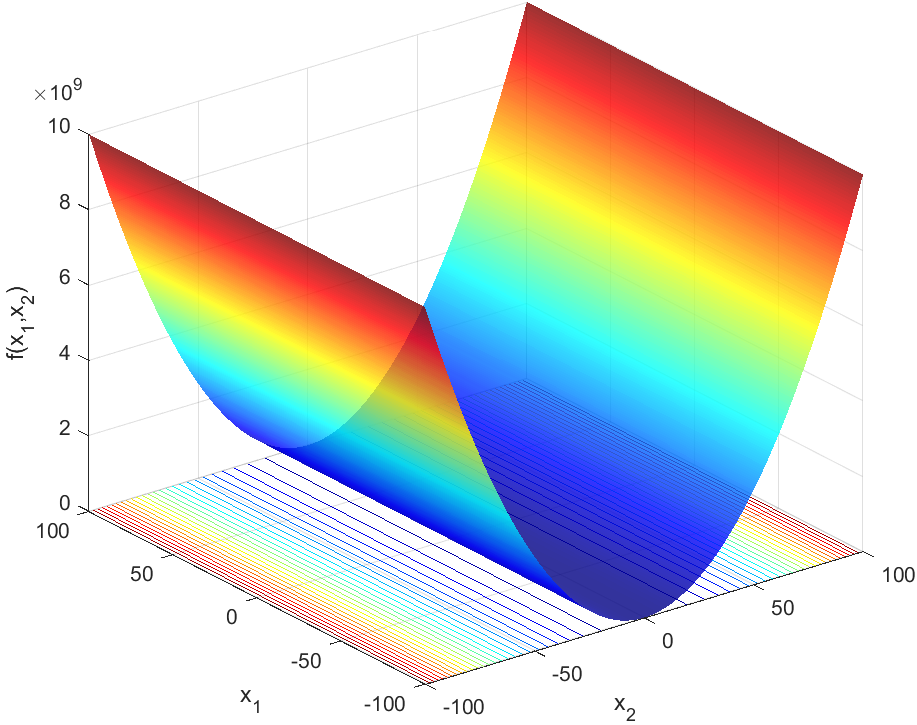}\label{fig:f3}}
\end{tabular}
\begin{tabular}{cc}
 \subfigure[{\scriptsize $f_4$.}]{\includegraphics[width=0.30\linewidth]{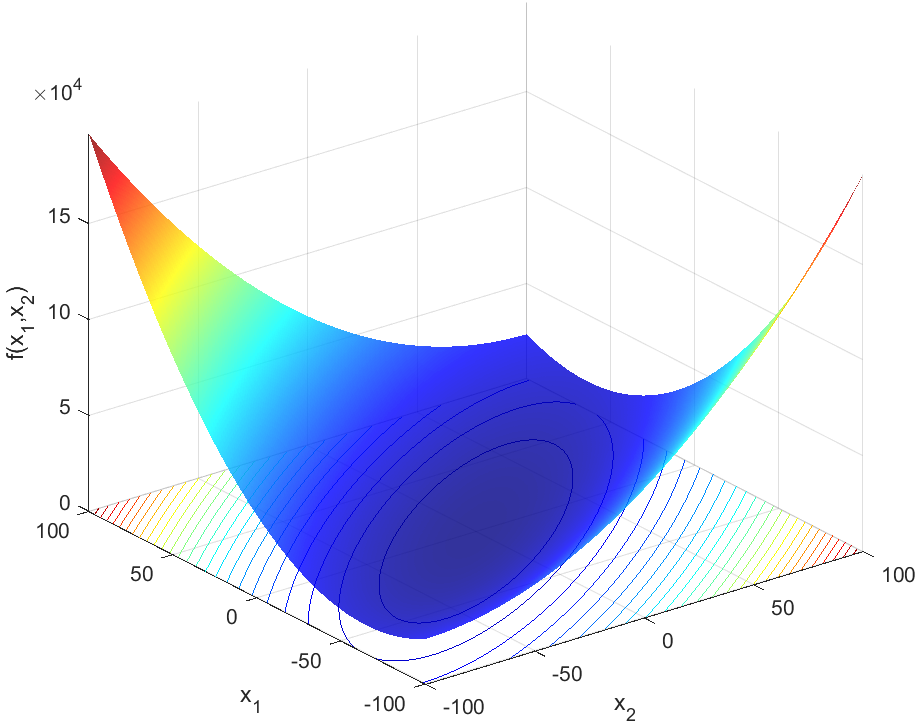}\label{fig:f4}}
&
\subfigure[{\scriptsize $f_5$ and $f_6$.}]{\includegraphics[width=0.30\linewidth]{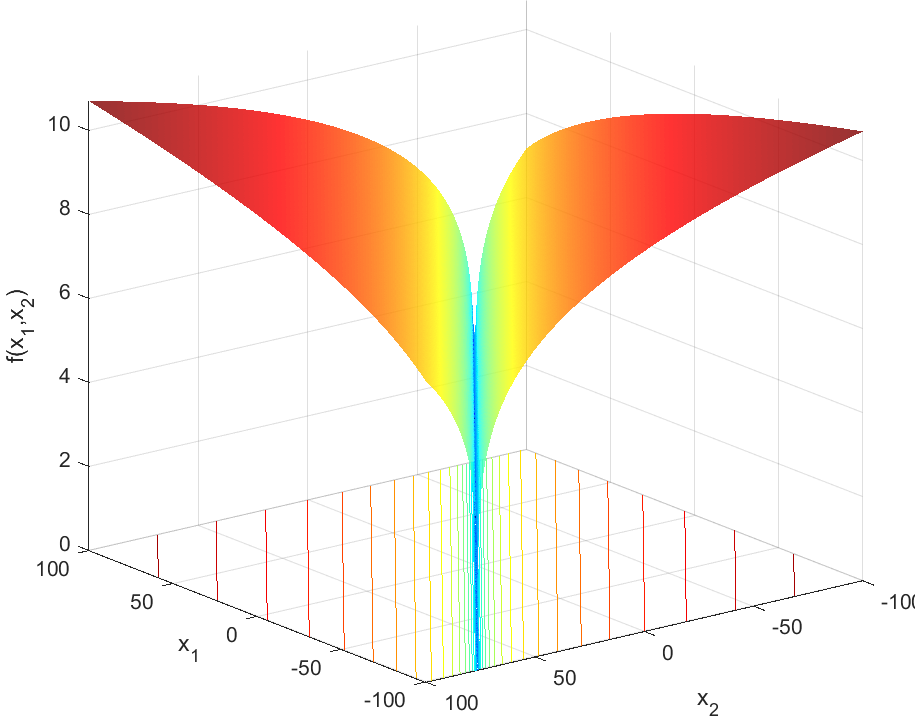}\label{fig:f5}}
\end{tabular}
\caption{Visualization of the 2-dimensional search space for unimodal functions $f_1$ to $f_6$. 
To enhance the clarity of the morphology depiction, we set $\sigma$ to zero and $\vec{m}$ to [0,0].}
\label{fig:Unimodal}
\end{figure*}

With $f_5$, our objective is to critically assess how optimization algorithms grapple with a triad of complexities. 
This problem instance stands as a challengng unimodal test, designed to spotlight algorithmic prowess (or lack thereof) in navigating intricate, narrow search topographies. Figure~\ref{fig:f5} offers a visual representation of $f_5$ in a 2-dimensional domain.

\subsubsection{$f_6$}

For $f_6$, we employ a configuration analogous to $f_5$. 
However, a key distinction lies in the variable interaction structure: $f_6$ boasts the maximum connectivity, rendering it a fully-connected, non-separable function. 
By juxtaposing the outcomes on $f_6$ with those on $f_5$, we can discern how algorithmic performance is modulated by the intricacies of connectivity structures in non-separable functions. Given that in a 2-dimensional domain, $f_6$ and $f_5$ exhibit marked similarities, we have opted to forgo a separate illustration for $f_6$.

\subsection{Multimodal Problem Instances with a Single Component}
\label{sec:sec:MultimodalScenarios}

This section introduces a series of multimodal problem instances, labeled as $f_7$ to $f_{15}$, each characterized by a single component in their landscape (i.e., $o=1$). 
Each instance is designed in a 30-dimensional space ($d=30$), with the minimum position, $\vec{m}$, chosen randomly within the range $[-80, 80]^d$. 
The global optimum value, $\sigma$, is randomly selected from the range $[-1,200,0]$. 
To ensure consistent and reproducible performance evaluations, the random seed employed for generating both the minimum position $\vec{m}$ and value $\sigma$ is kept constant across all instances. 
This consistency aids in offering reliable comparisons of optimization algorithms across these multimodal scenarios. 
In the subsequent sections, we delve into each instance, providing details on their unique attributes, variable interactions, and other relevant characteristics to enable a thorough analysis of optimization algorithms on these landscapes.

\subsubsection{$f_7$}

For $f_7$, we set the parameters as follows: $\lambda=1$, $\bm\mu=(0.2,0.2)$, $\bm\omega=(20,20,20,20)$, and $\mathbf{H}=\mathbf{R}=\mathbf{I}_{d \times d}$. 
This specific multimodal problem instance is distinguished by its symmetry, well-conditioning, and separability. 
Its separable characteristic enables optimization algorithms to exploit the problem's structure, optimizing each variable individually. 
A comparison of results from $f_7$ and $f_1$ can shed light on the challenges introduced by the presence of local optima in the basin of a component. 
Note that the ruggedness and modality of $f_7$ are low-moderate.
Figure~\ref{fig:f7} visualizes a 2-dimensional representation of $f_7$.

\subsubsection{$f_8$}

The configuration of $f_8$ is similar to that of $f_7$, except $\bm\omega=(50,50,50,50)$.
Similar to $f_7$, $f_8$ is multimodal, symmetric, well-conditioned, and separable.
However, in contrast to $f_7$, $f_8$ contains a significantly larger number of basin local optima, and these optima are narrower than those present in $f_7$.
By comparing the results obtained by the algorithms on $f_8$ with those on $f_7$, we can observe the impact of having a larger number of narrower local optima.
This characteristic poses additional challenges for some optimization algorithms, as they must navigate through a more intricate landscape with increased roughness, irregularities, and local optima. 
Figure~\ref{fig:f8} depicts a visual representation of $f_8$ in a 2-dimensional domain.

\subsubsection{$f_9$}
The configuration of $f_9$ is similar to that of $f_7$, but $\bm\mu=(1,1)$.
Similar to $f_7$, $f_9$ is also multimodal, symmetric, well-conditioned, and separable.
One notable difference in $f_9$ compared to $f_7$ is the significantly greater depth of its basin local optima.
By comparing the results obtained by the algorithms on $f_9$ with those on $f_7$, we can analyze the impact of the depth of the basin local optima on the algorithms' performance. 
The deeper optima in $f_9$ pose a greater challenge for some optimization algorithms, as they need to exploit deeper in the local optima, which may result in loss of diversity and trapping in them. 
Figure~\ref{fig:f9} visualizes a 2-dimensional representation of $f_9$.

\subsubsection{$f_{10}$}

For $f_{10}$, we aim to examine the influence of an asymmetric basin of attraction. 
To this end, we configure this instance with parameters $\bm\mu=(0.2,0.5)$ and $\bm\omega=(20,50,10,25)$, while retaining the other parameters consistent with $f_7$ to $f_9$.  
The specific values chosen for $\bm\mu$ and $\bm\omega$ lead to a unique feature: when observing the basin of $f_{10}$ from various directions, the intricacies introduced by local optima lie somewhere between those observed in $f_{7}$ to $f_{9}$. 
The asymmetric nature of the basin, sets $f_{10}$ apart. 
By juxtaposing algorithmic results on $f_{10}$ against those from $f_{7}$ to $f_{9}$, we can discern the challenges posed by an asymmetric basin of attraction and its nuanced interplay with local optima.
Figure~\ref{fig:f10} shows a visual representation of $f_{10}$ in a 2-dimensional domain.

\subsubsection{$f_{11}$}

$f_{11}$ is configured as a non-separable problem instance. 
While it shares similarities with $f_{10}$, a key distinction is its fully connected variable interaction structure. 
For $f_{11}$, the rotation matrix $\mathbf{R}$ is derived using Algorithm~\ref{alg:RotationControlled}, where $\mathfrak{p}=1$ and every angle in $\bm\Theta$ is uniformly and randomly selected from the range $(- \pi, \pi)$. 
A comparison between the algorithmic outcomes on $f_{10}$ and $f_{11}$ elucidates the impact of transitioning the variable interaction structure from separable to fully non-separable. 
Specifically, it illustrates how an algorithm navigates an asymmetric multimodal space where each variable interacts with every other variable.
Figure~\ref{fig:f11} illustrates a 2-dimensional projection of $f_{11}$.

\subsubsection{$f_{12}$}

This problem instance is a partially separable version of $f_{10}$, and exhibits a unique variable interaction structure, consisting of three groups of variables with fully connected interactions within each group. 
Each group contains 10 randomly selected variables from the 30 variables of the problem instance.
Notably, no variable from one group interacts with any variable from the other two groups. 
For each group, we use a distinct angle from $\{\frac{\pi}{4}, \frac{3\pi}{4}, \frac{\pi}{8}\}$ in $\bm\Theta$ to establish variable interactions.
The purpose of introducing $f_{12}$ is to evaluate whether the algorithms can effectively identify and exploit the partial separable variable interaction structure in the problem instances.
By comparing the results obtained on $f_{12}$ and $f_{11}$ we gain insights into how well the algorithms can adapt to and benefit from the partially separable variable interaction of $f_{12}$ in comparison to the non-separable $f_{11}$.
Since in a 2-dimensional domain, $f_{12}$ and $f_{11}$ are similar, we do not provide a separate figure for $f_{12}$.

\subsubsection{$f_{13}$}

This problem instance is designed to critically test the exploration capability of optimization algorithms within a non-separable setting. Its configuration parameters include: $\lambda=1$, $\bm\mu=(1,1)$, $\bm\omega=(50,50,50,50)$, and $\mathbf{H}=\mathbf{I}_{d \times d}$.
The rotation matrix, $\mathbf{R}$, is formulated using Algorithm~\ref{alg:RotationControlled}, with $\mathfrak{p}=1$. 
Furthermore, each angle within $\bm\Theta$ is uniformly and randomly determined from the interval $(- \pi, \pi)$. 
Characteristically, $f_{13}$ is well-conditioned, symmetric, and highly rugged. 
Its fully connected non-separable nature is punctuated by numerous profound local optima, intensifying the optimization challenge.
Figure~\ref{fig:f13} shows a projection of $f_{13}$ in a 2-dimensional space.

\subsubsection{$f_{14}$}

This problem instance is designed to encompass several challenging problem characteristics simultaneously.
The configuration of this problem instance is as follows: $\lambda$ is set to 0.6, $\bm\mu=(0.7,0.2)$, $\bm\omega=(25,10,20,50)$.
the rotation matrix $\mathbf{R}$ is derived using Algorithm~\ref{alg:RotationControlled}, where $\mathfrak{p}=1$ and every angle in $\bm\Theta$ is uniformly and randomly selected from the range $(- \pi, \pi)$.
Additionally, we set two randomly chosen elements of the principal diagonal of $\mathbf{H}$ to 0.01 and $[10^3]$, while the remaining elements are randomly chosen (with uniform distribution) from the range $[1,10^3]$.
The characteristics of $f_{14}$ include being ill-conditioned (condition number of $\mathbf{H}$ is $[10^5]$, but the condition number of the function id lower as it is damped by the value of $\lambda$), super-linear, asymmetric, rough, irregular, and non-separable with a fully connected variable interaction structure.
Figure~\ref{fig:f14} illustrates a 2-dimensional representation of $f_{14}$.

\subsubsection{$f_{15}$}

$f_{15}$ is crafted to represent a challenging optimization landscape characterized by multiple valleys, with the optimal solution ensconced within a slender valley that contains local minima. 
This configuration intensifies the optimization challenge, as algorithms must adeptly navigate the rugged and convoluted terrain to pinpoint the global optimum. 
Specifically, the configuration parameters for this instance are set as: $\lambda=0.1$, $\bm\mu=(1,1)$, and $\bm\omega=(10,10,10,10)$. The rotation matrix $\mathbf{R}$ is computed using Algorithm~\ref{alg:RotationControlled}, with $\mathfrak{p}=1$. 
Each angle in $\bm\Theta$ is uniformly and randomly determined from the interval $(- \pi, \pi)$. 

In addition, two elements from the principal diagonal of $\mathbf{H}$ are set to values of 1 and $10^5$ respectively, while the remaining elements are uniformly selected from the range $[1,10^5]$.
 
\begin{figure*}[!t]
\centering
\begin{tabular}{ccc}
    \subfigure[{\scriptsize $f_7$}]{\includegraphics[width=0.30\linewidth]{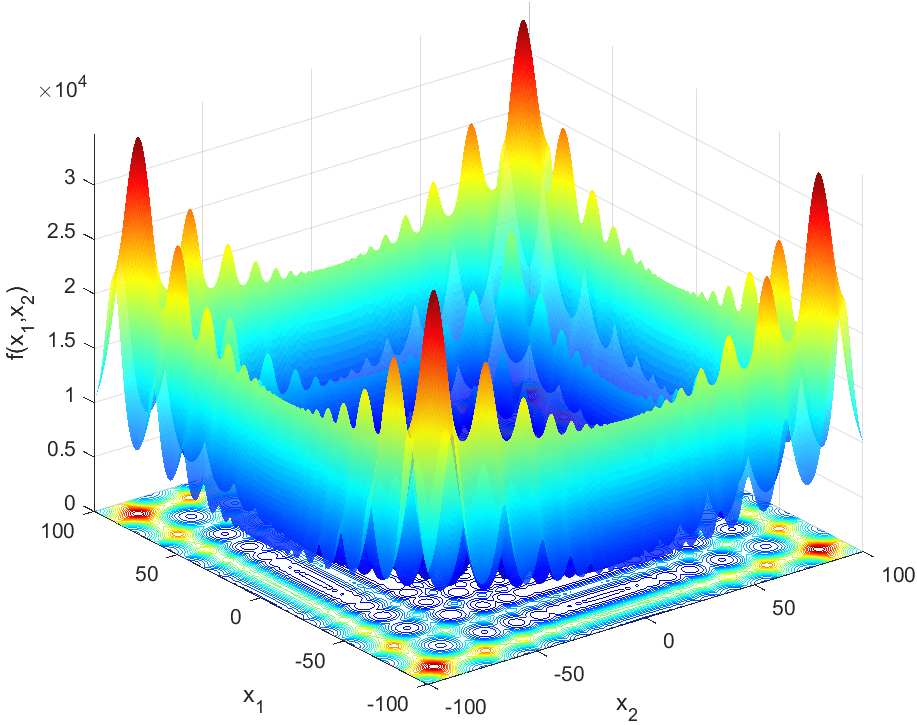}\label{fig:f7}}
&
     \subfigure[{\scriptsize $f_8$}]{\includegraphics[width=0.30\linewidth]{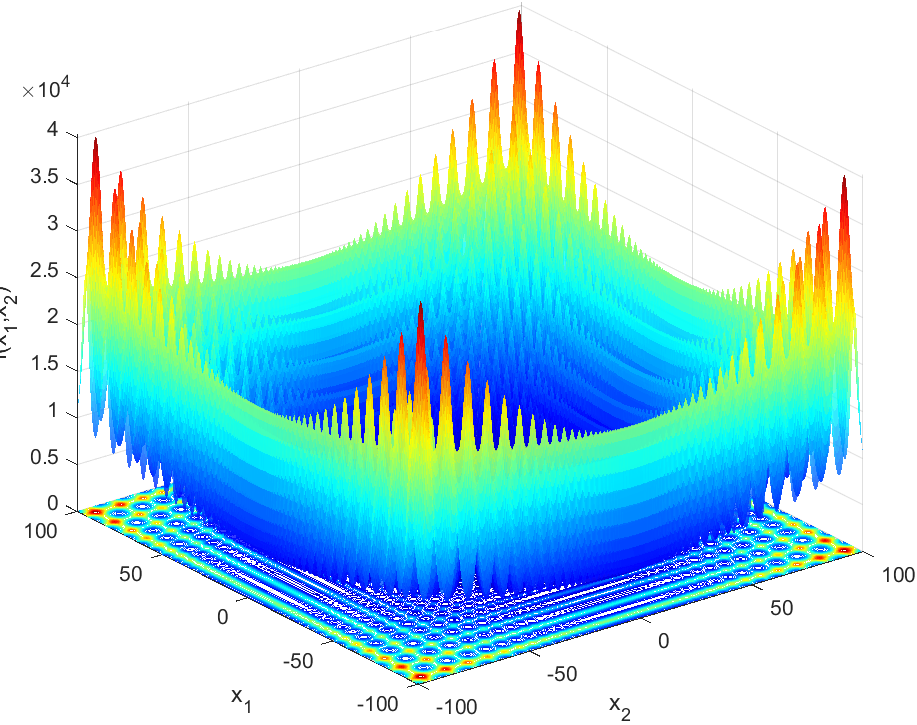}\label{fig:f8}}
&
\subfigure[{\scriptsize $f_9$}]{\includegraphics[width=0.30\linewidth]{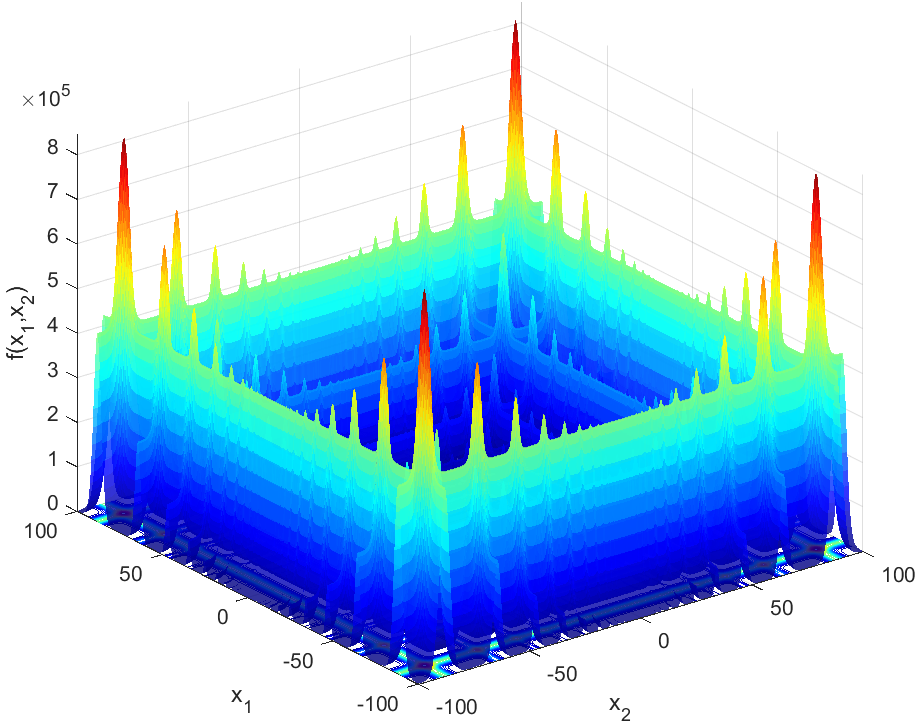}\label{fig:f9}}
\end{tabular}
\begin{tabular}{ccc}
    \subfigure[{\scriptsize $f_{10}$}]{\includegraphics[width=0.30\linewidth]{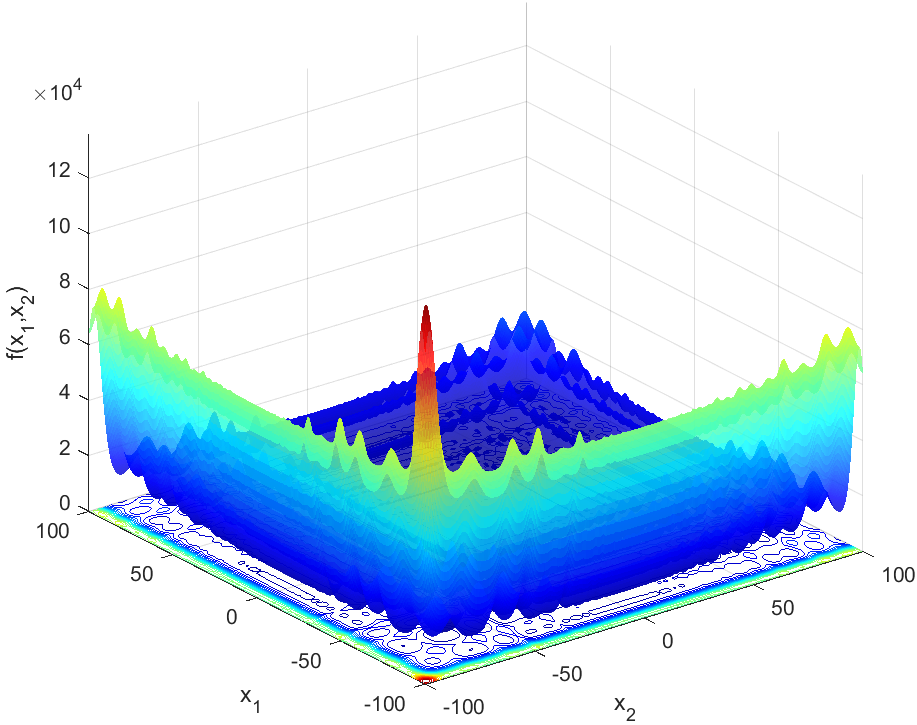}\label{fig:f10}}
&
     \subfigure[{\scriptsize $f_{11}$ and $f_{12}$}]{\includegraphics[width=0.30\linewidth]{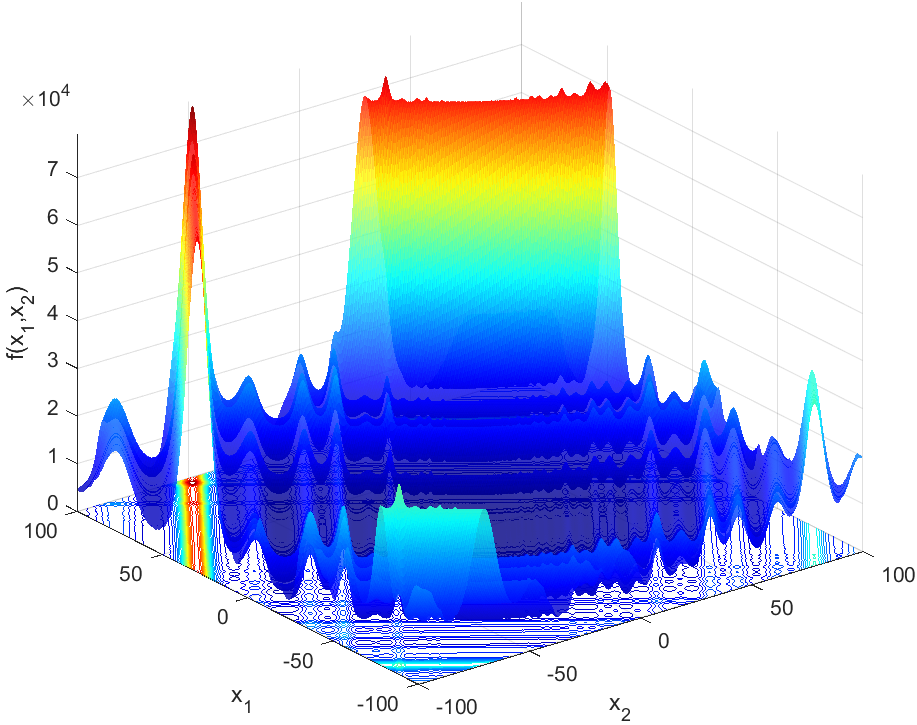}\label{fig:f11}}
&
\subfigure[{\scriptsize $f_{13}$}]{\includegraphics[width=0.30\linewidth]{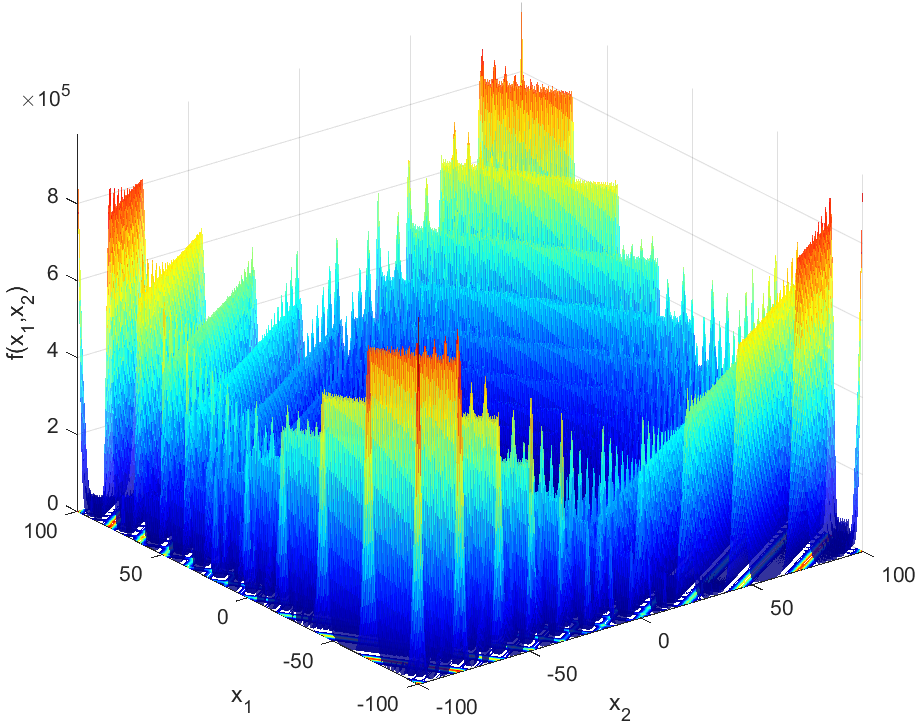}\label{fig:f13}}
\end{tabular}
\begin{tabular}{cc}
 \subfigure[{\scriptsize $f_{14}$}]{\includegraphics[width=0.30\linewidth]{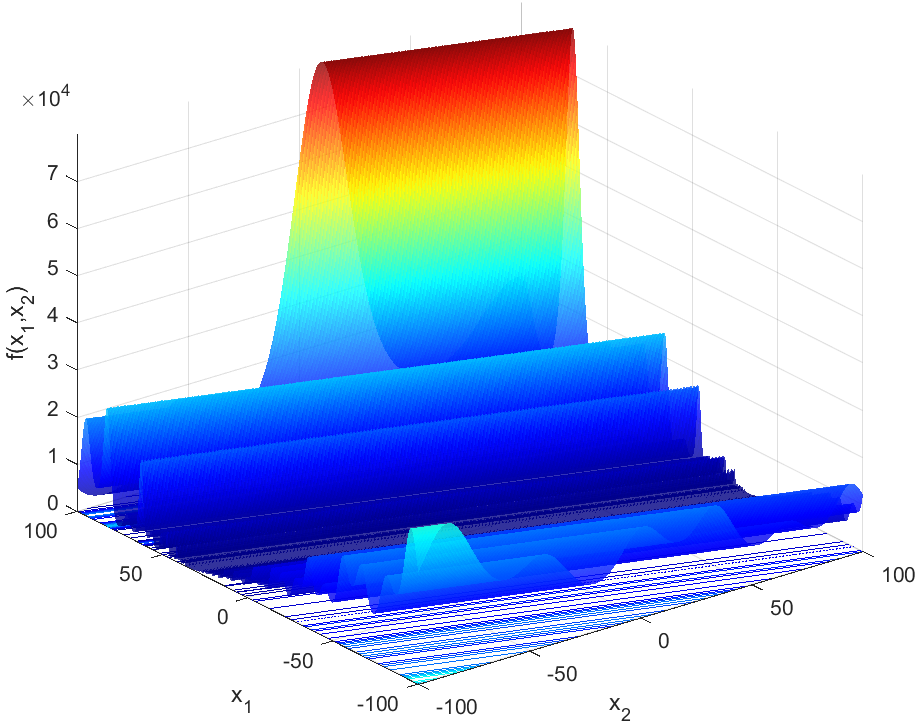}\label{fig:f14}}
&
\subfigure[{\scriptsize $f_{15}$}]{\includegraphics[width=0.30\linewidth]{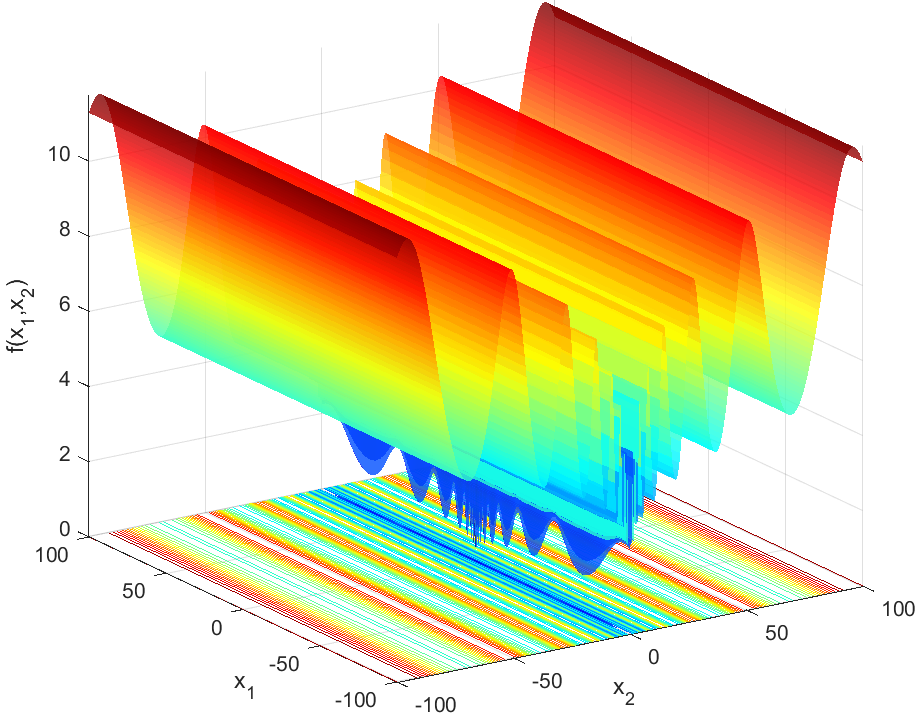}\label{fig:f15}}
\end{tabular}
\caption{Visualization of the 2-dimensional search space for multimodal functions $f_7$ to $f_{15}$. 
To enhance the clarity of the morphology depiction, we set $\sigma$ to zero and $\vec{m}$ to [0,0].}
\label{fig:Multimodal}
\end{figure*}

In addition, two elements from the principal diagonal of $\mathbf{H}$ are set to values of 1 and $10^5$ respectively. The remaining elements are then selected from a heavy-tail distribution by using a beta distribution with both $\alpha$ and $\beta$ parameters set to 0.2, spanning the range $[1,10^5]$.
Notably, $f_{15}$ exhibits properties of ill-conditioning, pronounced sub-linearity, asymmetry, and full-connectivity in its non-separable variable interaction structure.
Figure~\ref{fig:f15} depicts a 2-dimensional projection of $f_{15}$.

\subsection{Multimodal Problem Instances with Multiple Components}
\label{sec:sec:MultiOptimaScenarios}

This section introduces a suite of problem instances characterized by landscapes containing multiple components (i.e., $o \geq 2$).
These landscapes present intricate terrains populated by multiple competing optima, each distinguished by its unique morphological attributes and basin of attraction. 
This setup poses a substantial exploration challenge for optimization algorithms. 
Notably, each component has an associated $\sigma_i \leq 0$, with the global optimum value being $\min_{i\in\{1,\ldots,o\}}(\sigma_i)$, and the global optimum position corresponding to the minimum position ($\vec{m}$) of the component with the smallest $\sigma$. 
It is essential to highlight that all problem instances in this category inherently manifest non-separability due to the coexistence of multiple components. 
This non-separability persists irrespective of the particular variable interaction configurations within individual components~\cite{yazdani2018thesis,yazdani2019scaling}.

\subsubsection{$f_{16}$}

For $f_{16}$, we design a landscape encompassing five homogeneous components. 
Each component is characterized by unimodality, smoothness, regularity, symmetry, full separability (despite the overall landscape's non-separability), and well-conditioning.
For each component $i \in \{1, \ldots, 5\}$, the configuration is delineated as follows: $\bm\mu_i=(0,0)$, $\bm\omega_i=(0,0,0,0)$, and both $\mathbf{R}_i$ and $\mathbf{H}$ are set to $\mathbf{I}_{d \times d}$, with $\lambda_i=1$. 
One component has its $\sigma$ value set to -5,000, while the $\sigma$ values for the remaining components are randomly determined via $\mathcal{U}(-4,500,-4,000)$, where $\mathcal{U}(a,b)$ uniformly generates a random number within the interval $(a,b)$. 
The minimum position, $\vec{m}_i$, for each component is randomly drawn from the range $[-80,80]^{d}$.

Importantly, all parameter settings for each component represent the most straightforward configurations, purposefully chosen. 
This deliberate design aims to evaluate the algorithms' adeptness in navigating the landscape, gauging their resistance against being captivated by the allure of expansive and high-quality competing components.
Figure~\ref{fig:f16} provides an illustrative example of $f_{16}$ in a 2-dimensional domain.

\subsubsection{$f_{17}$}

In $f_{17}$, we build upon the foundational structure of $f_{16}$, incorporating added intricacies through ill-conditioning and non-separability. 
For each component $i$, elements of the principal diagonal of $\mathbf{H}_i$ are randomly chosen using $\mathcal{U}(0.01,100)$. 
By setting $\mathfrak{p}_i=0.5$, we ensure that each component has different variable interaction structure and connectivity. 
Additionally, each angle within $\bm\Theta_i$ is determined using $\mathcal{U}(-\pi,\pi)$.
As a result, the basin of attraction for each component in $f_{17}$ exhibits varied condition numbers and distinct variable interactions, giving rise to a landscape marked by its heterogeneous components. 
Such heterogeneity necessitates that, when transitioning between basins in $f_{17}$, optimization algorithms adjust both their convergence direction and step size across different dimensions.

While the configurations for $\mathbf{H}_i$ and $\mathbf{R}_i$ distinguish $f_{17}$ from $f_{16}$, the remaining parameters align with those of $f_{16}$. 
The primary objective of introducing $f_{17}$ is to assess the resilience and adaptability of optimization algorithms when confronted with varied condition numbers and intricate variable interaction structures.
An example landscape of $f_{17}$ is visually depicted in Figure~\ref{fig:f17}, demonstrating characteristics in a 2-dimensional space.

\subsubsection{$f_{18}$}

This problem instance is configured similarly to $f_{16}$, but with distinct variations. 
For each component, the values for $\mu_{1,i}$ and $\mu_{2,i}$ are chosen using $\mathcal{U}(0.2,0.5)$, and values for $\omega_{1,i}$ through $\omega_{4,i}$ are randomly chosen by $\mathcal{U}(5,50)$. 
We set $\mathfrak{p}_i=0.5$, therefore, each component has different variable interaction structure and connectivity. 
In addition, each angle within $\bm\Theta_i$ is determined using $\mathcal{U}(-\pi,\pi)$.

Consequently, $f_{18}$ showcases a highly multimodal nature with rough, irregular landscapes, along with pronounced local and global asymmetry. 
While each component is fully separable, they exhibit pronounced differences in the morphologies of their local optima basins and patterns. Such characteristics further challenge optimization algorithms to navigate through intricate terrains laden with a multitude of local optima, each possessing distinct shapes and basin attributes.
Figure~\ref{fig:f18} showcases a representative 2-dimensional visualization of $f_{18}$.

\subsubsection{$f_{19}$}

In $f_{19}$, we maintain a configuration similar to $f_{18}$, but with specific alterations. 
For each component, both $\mu_{1,i}$ and $\mu_{2,i}$ are set to 0.5, and each $\omega_{1,i}$ to $\omega_{4,i}$ value is chosen by $\mathcal{U}(50,100)$. 
This particular setting of the $\bm\mu$ and $\bm\omega$ parameters leads to the creation of an extremely rough and rugged landscape with enormous number of local optima. 
A sample 2-dimensional representation of $f_{19}$ is presented in Figure~\ref{fig:f19}, offering insights into its  morphological characteristics.

\subsubsection{$f_{20}$}

In $f_{20}$, we retain a configuration similar to $f_{18}$, with the following modifications: $\lambda_i$ for each component is set to 0.25, the minimum position $\vec{m}_i$ for each component is randomly chosen in a small sub region limited in $[-75, -25]^{d}$, and for one component $\sigma$ value is set to -100, while the $\sigma$ values for the remaining components are randomly determined via $\mathcal{U}(-99,-98)$. 

Owing to these adjustments in parameter settings, $f_{20}$ displays components with a highly sub-linear basin. 
Moreover, all the components in $f_{20}$ are closely packed together within a small sub-region, which occupies $\frac{1}{4^d}$ of the entire search space.
The presence of multiple sharp components in a confined space challenges certain algorithms' ability to transition between components, especially after they have converged in this small sub-region and lost their diversity. 
Furthermore, $f_{20}$ showcases heterogeneous components due to the varying basin local optima patterns and morphologies associated with each. 
This problem instance serves as a distinctive test case for evaluating the algorithms' efficacy in situations where several competing components are in close quarters.
In Figure~\ref{fig:f20}, a 2-dimensional example of $f_{20}$ is depicted.

\subsubsection{$f_{21}$}

In this problem instance, we present a landscape that comprises five heterogeneous components. 
For each component $i\in\{1,\ldots,5\}$, the configuration is set as: $\lambda_i=0.5$, the values for $\mu_{1,i}$ and $\mu_{2,i}$ are chosen using $\mathcal{U}(0.1,0.2)$, and values for $\omega_{1,i}$ through $\omega_{4,i}$ are randomly chosen by $\mathcal{U}(5,10)$. 
We set $\mathfrak{p}_i=0.5$, therefore, each component has different variable interaction structure and connectivity. 
In addition, each angle within $\bm\Theta_i$ is determined using $\mathcal{U}(-\pi,\pi)$.
The $\sigma$ values for the components are set to $[-50,-45,-40,-40,-40]$ respectively. 
The second component (with $\sigma=-45$) is centrally positioned in the search space, i.e., $\vec{m}_2 = (0,0,\ldots,0)$, and has $\mathbf{H}_2=\mathrm{diag}(1,1,\ldots,1)$. 
For the other components, $\mathbf{H}_i=\mathrm{diag}(5,5,\ldots,5)$ and their minimum positions are randomly selected to lie outside the domain $[-30,30]^{d}$ but within $[-90,90]^{d}$.

The landscape of $f_{21}$ showcases a heterogeneity in the sizes of the basins of attraction among its components. Notably, a broad promising region situated at the center of the search space spans approximately half of it. This dominant basin is around five times larger than the others, making $f_{21}$ a deceptive problem. Such a landscape challenges optimization algorithms to bypass the extensive central basin and discover the other narrower components. Given its deceptive nature, $f_{21}$ serves as a robust test case for evaluating algorithmic performance in navigating complex landscapes with heterogeneous basins of attraction.
An illustrative 2-dimensional example of $f_{21}$ is given in Figure~\ref{fig:f21}.

\subsubsection{$f_{22}$}

In $f_{22}$, the landscape features two distant multimodal components located at opposite ends of the search space. 
For the first component, $\vec{m}_1$ is chosen using $\mathcal{U}(80,90)^d$, while for the second component, $\vec{m}_2$ is chosen using $\mathcal{U}(-80,-90)^d$. 
$\sigma_1$ and $\sigma_2$ are  set to -1,000 and -950, respectively. 
For both components, the $\lambda_1$ and $\lambda_2$ values are set to 1 and 0.9 respectively. 
Additionally, each principal diagonal element of $\mathbf{H}_1$ and $\mathbf{H}_2$ is chosen using $\mathcal{U}(1,10)$. 
Moreover, $\mathfrak{p}_1 = \mathfrak{p}_2 = 0.7$ and the angle values are chosen using $\mathcal{U}(-\pi, \pi)$. 
Furthermore, the $\bm\mu_1$ and $\bm\mu_2$ values are set to [0.5,0.5], and the $\bm\omega$ values are determined by $\mathcal{U}(20,50)$. 
$f_{22}$ represents a deceptive problem instance highlighted by the presence of two high-quality, distant components. 
Using a smaller $\lambda$ value for the inferior component expands its basin of attraction, which in turn increases the deceptiveness of this instance.
In Figure~\ref{fig:f22}, a 2-dimensional landscape generated by $f_{22}$ is depicted.

\subsubsection{$f_{23}$}

In $f_{23}$, the landscape is generated by composing five overlapped components. 
Each component has $\lambda_i$ set to 0.4, and $\mathbf{H}_i = \mathbf{I}_{d \times d}$. 
The minimum position $\vec{m}_i$ for all components is identical and is randomly chosen using $\mathcal{U}(-80,80)^d$. 
Additionally, the $\sigma_i$ values for all components are set to -100. 
For all components, $\mathfrak{p}_{i} = 0.75$ and the angle values for each $\bm\Theta_i$ are chosen randomly using $\mathcal{U}(-\pi, \pi)$, ensuring that each component is uniquely rotated in different planes. 
Moreover, for each component, both $\mu_{1,i}$ and $\mu_{2,i}$ are set to 0.5, and each of $\omega_{1,i}$ to $\omega_{4,i}$ is randomly chosen using $\mathcal{U}(20,50)$.
These five overlapped and randomly rotated components yield a single, highly irregular and rough visible component in the landscape. 
This results in a complex local optima pattern, wherein the variable interaction structure across these local optima is diverse. 
A 2-dimensional landscape generated by $f_{23}$ is illustrated in Figure~\ref{fig:f23}.

\subsubsection{$f_{24}$}

\begin{figure*}[!t]
\centering
\begin{tabular}{ccc}
    \subfigure[{\scriptsize $f_{16}$}]{\includegraphics[width=0.30\linewidth]{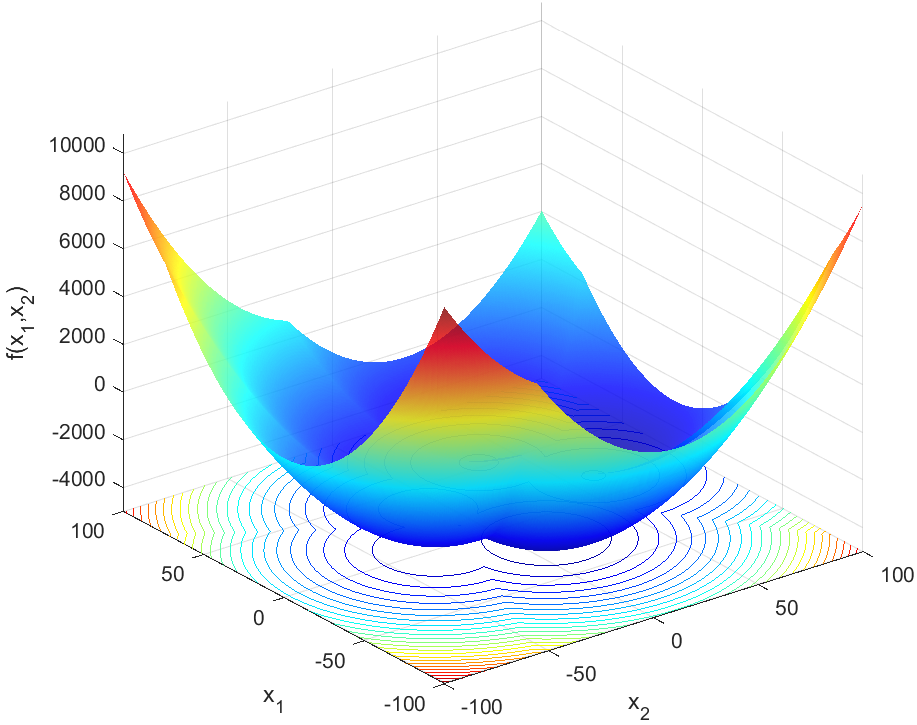}\label{fig:f16}}
&
     \subfigure[{\scriptsize $f_{17}$}]{\includegraphics[width=0.30\linewidth]{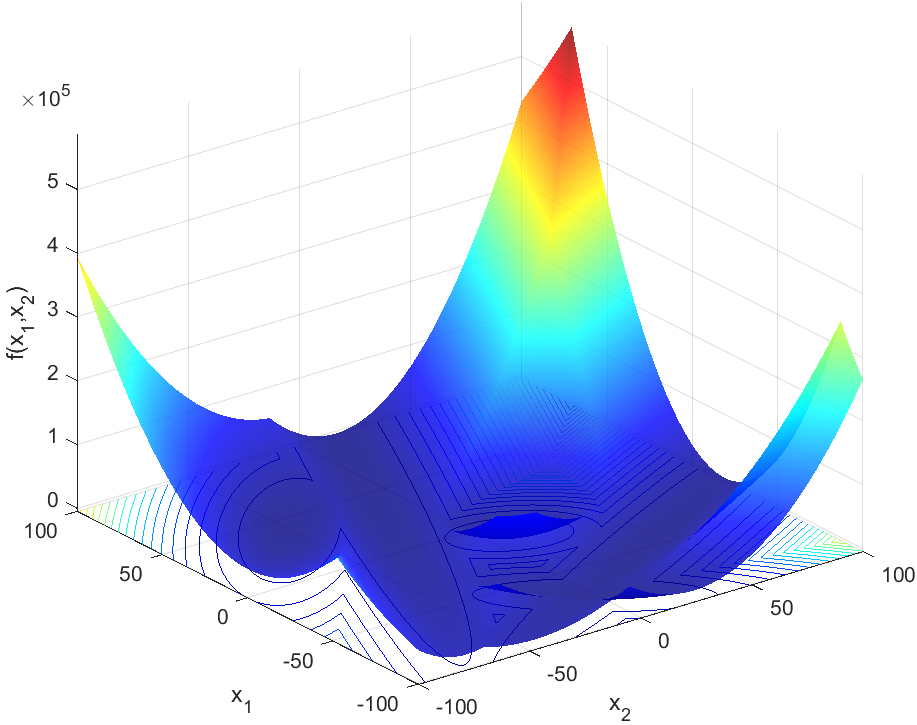}\label{fig:f17}}
&
\subfigure[{\scriptsize $f_{18}$}]{\includegraphics[width=0.30\linewidth]{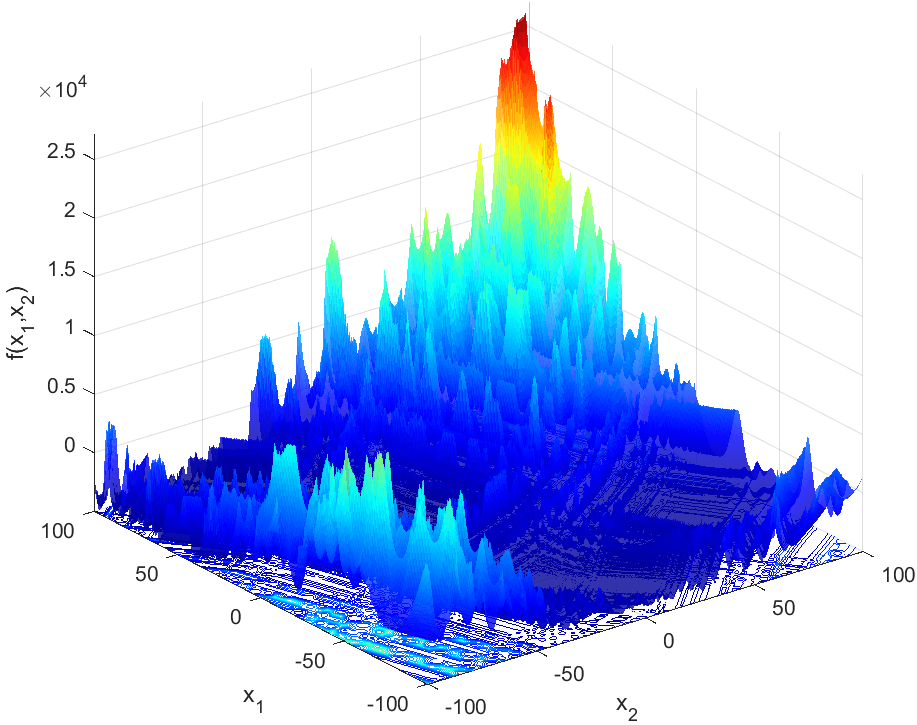}\label{fig:f18}}
\\
    \subfigure[{\scriptsize $f_{19}$}]{\includegraphics[width=0.30\linewidth]{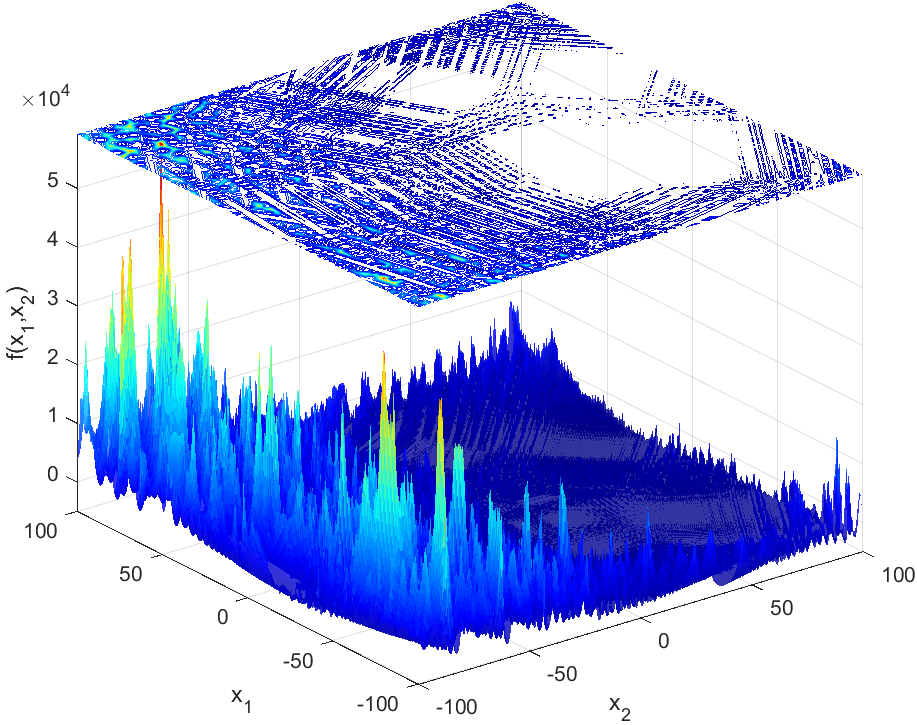}\label{fig:f19}}
&
     \subfigure[{\scriptsize $f_{20}$}]{\includegraphics[width=0.30\linewidth]{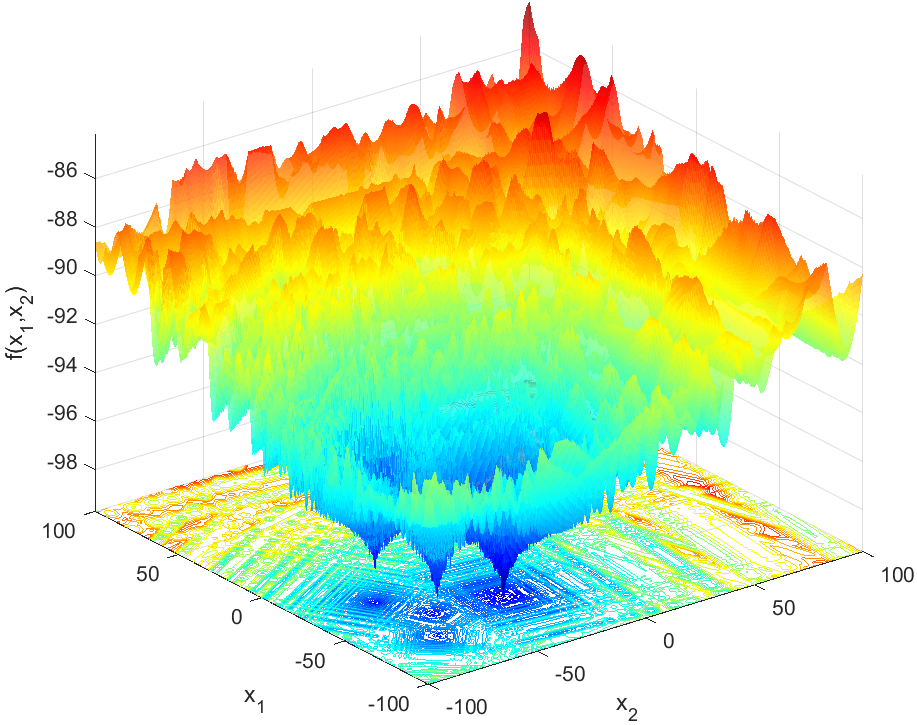}\label{fig:f20}}
&
\subfigure[{\scriptsize $f_{21}$}]{\includegraphics[width=0.30\linewidth]{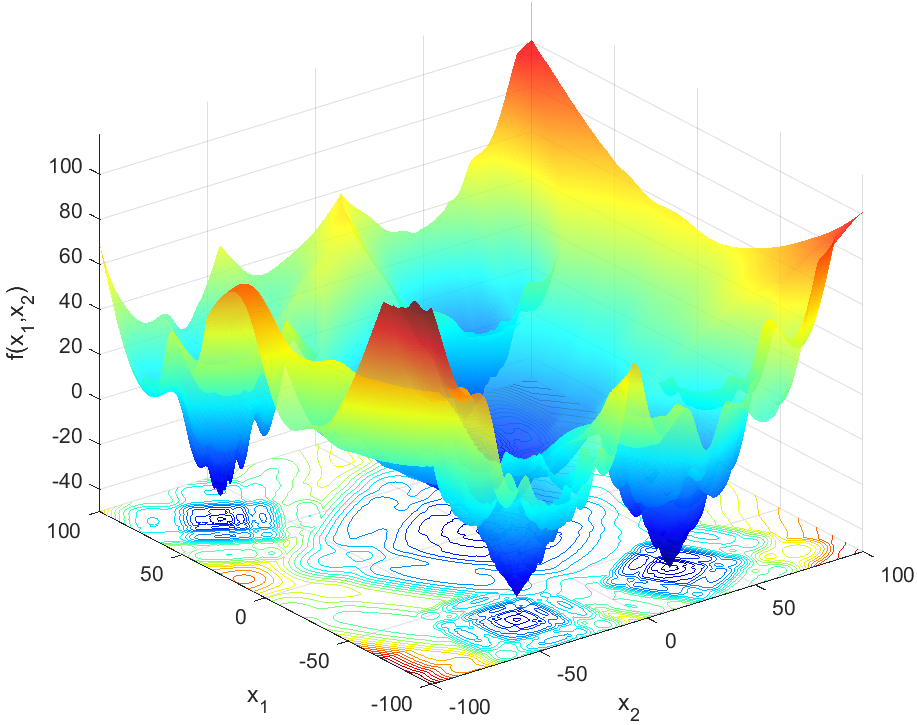}\label{fig:f21}}
\\
    \subfigure[{\scriptsize $f_{22}$}]{\includegraphics[width=0.30\linewidth]{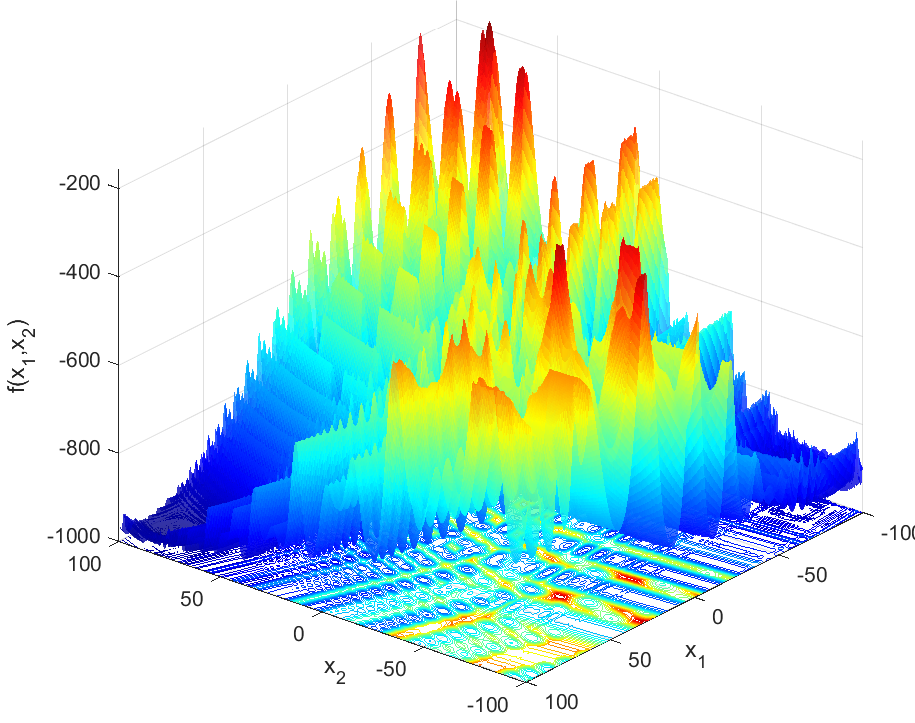}\label{fig:f22}}
&
     \subfigure[{\scriptsize $f_{23}$}]{\includegraphics[width=0.30\linewidth]{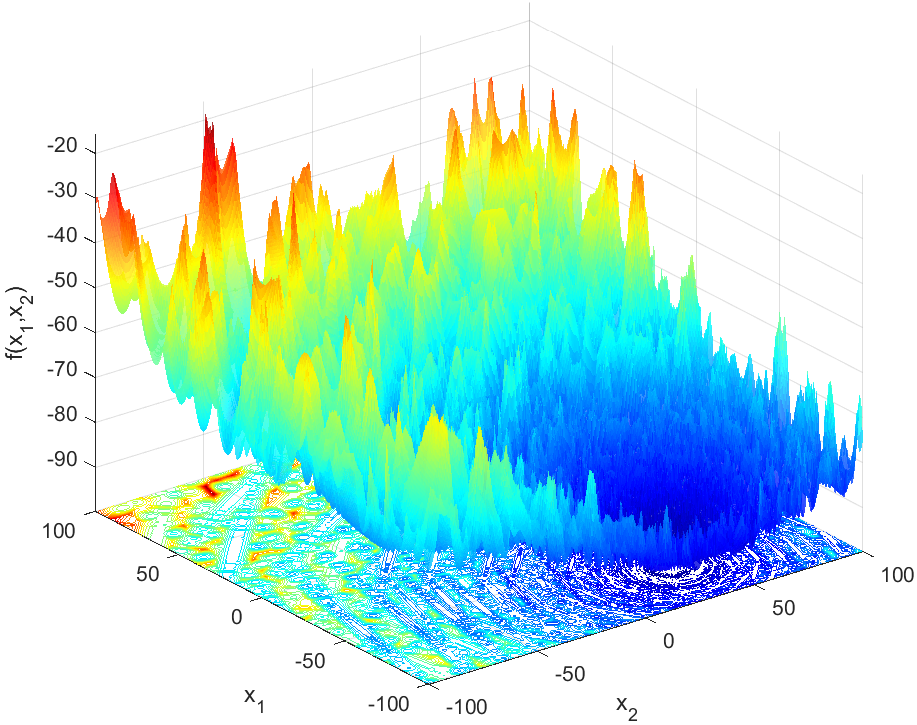}\label{fig:f23}}
&
\subfigure[{\scriptsize $f_{24}$}]{\includegraphics[width=0.30\linewidth]{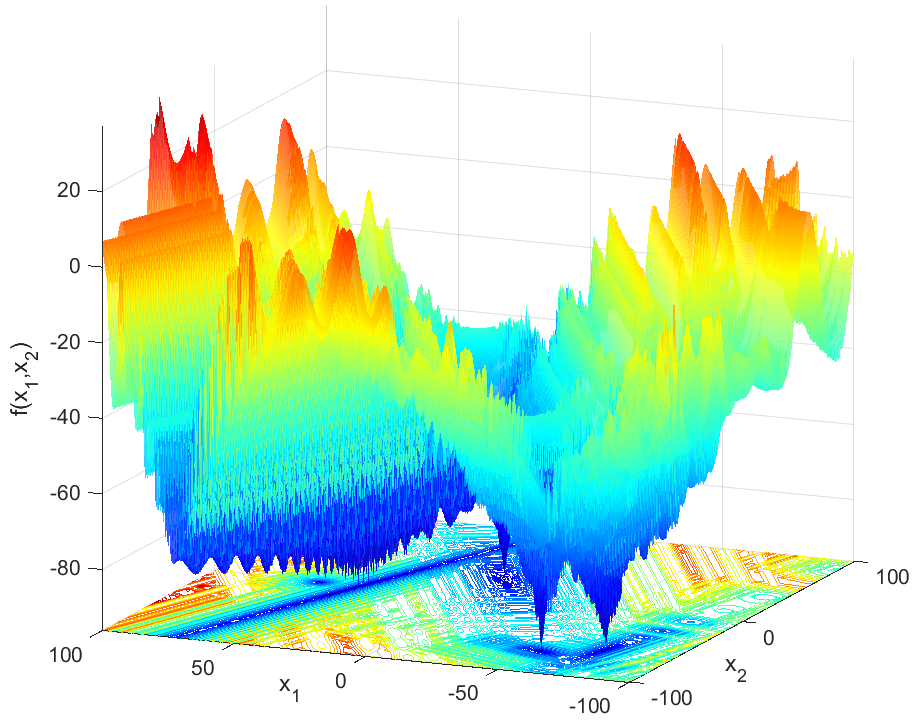}\label{fig:f24}}
\end{tabular}
\caption{Visualization of the 2-dimensional search space for multimodal functions with multiple components $f_{16}$ to $f_{24}$. }
\label{fig:MultiComponent}
\end{figure*}

This problem instance contains five heterogeneous components with various challenging characteristics. 
For each component, $\lambda_i$ is set to 0.25. 
One component is assigned a $\sigma$ of -100, while the $\sigma$ values for the remaining components are randomly determined via $\mathcal{U}(-99,-98)$.  
The minimum position, $\vec{m}_i$, for each component is randomly selected from the range $[-80,80]^d$.
The values for $\mu_{1,i}$ and $\mu_{2,i}$ are chosen using $\mathcal{U}(0.2,0.5)$, and values for $\omega_{1,i}$ through $\omega_{4,i}$ are randomly selected by $\mathcal{U}(5,50)$. 
Therefore, each component exhibits an asymmetric, rugged, and sub-linear basin. 
Moreover, for each component $i$, the elements of the principal diagonal of $\mathbf{H}_i$ are randomly chosen using $\mathcal{U}(1,10^5)$, leading to the generation of ill-conditioned components with different conditioning degrees. 
By setting $\mathfrak{p}_i = 0.75$, we ensure that each component exhibits a different non-separable variable interaction structure with a high degree of connectivity. 
Additionally, each angle within $\bm\Theta_i$ is determined using $\mathcal{U}(-\pi,\pi)$.
Due to these intricate characteristics, $f_{24}$ presents a highly challenging landscape where both the exploration and exploitation capabilities of optimization algorithms are rigorously tested.
Figure~\ref{fig:f24} shows a 2-dimensional landscape generated by $f_{24}$.

\small
\bibliography{bib} 

\begin{thebibliography}{10}
\providecommand{\url}[1]{#1}
\csname url@samestyle\endcsname
\providecommand{\newblock}{\relax}
\providecommand{\bibinfo}[2]{#2}
\providecommand{\BIBentrySTDinterwordspacing}{\spaceskip=0pt\relax}
\providecommand{\BIBentryALTinterwordstretchfactor}{4}
\providecommand{\BIBentryALTinterwordspacing}{\spaceskip=\fontdimen2\font plus
\BIBentryALTinterwordstretchfactor\fontdimen3\font minus
  \fontdimen4\font\relax}
\providecommand{\BIBforeignlanguage}[2]{{%
\expandafter\ifx\csname l@#1\endcsname\relax
\typeout{** WARNING: IEEEtran.bst: No hyphenation pattern has been}%
\typeout{** loaded for the language `#1'. Using the pattern for}%
\typeout{** the default language instead.}%
\else
\language=\csname l@#1\endcsname
\fi
#2}}
\providecommand{\BIBdecl}{\relax}
\BIBdecl

\bibitem{yazdani2023GNBGinstances}
\BIBentryALTinterwordspacing
D.~Yazdani, \emph{Generalized Numerical Benchmark Generator (GNBG)-Problem
  Instances (MATLAB Source Code)}, 2023 (accessed December 10, 2023). [Online].
  Available: \url{https://github.com/Danial-Yazdani/GNBG-Instances}
\BIBentrySTDinterwordspacing

\bibitem{beiranvand2017best}
V.~Beiranvand, W.~Hare, and Y.~Lucet, ``Best practices for comparing
  optimization algorithms,'' \emph{Optimization and Engineering}, vol.~18, pp.
  815--848, 2017.

\bibitem{hansen2021coco}
N.~Hansen, A.~Auger, R.~Ros, O.~Mersmann, T.~Tu{\v{s}}ar, and D.~Brockhoff,
  ``{COCO}: A platform for comparing continuous optimizers in a black-box
  setting,'' \emph{Optimization Methods and Software}, vol.~36, no.~1, pp.
  114--144, 2021.

\bibitem{mei2021structural}
L.~Mei and Q.~Wang, ``Structural optimization in civil engineering: a
  literature review,'' \emph{Buildings}, vol.~11, no.~2, p.~66, 2021.

\bibitem{yazdani2021DOPsurveyPartA}
D.~Yazdani, R.~Cheng, D.~Yazdani, J.~Branke, Y.~Jin, and X.~Yao, ``A survey of
  evolutionary continuous dynamic optimization over two decades -- part {A},''
  \emph{IEEE Transactions on Evolutionary Computation}, vol.~25, no.~4, pp.
  609--629, 2021.

\bibitem{yazdani2021DOPsurveyPartB}
------, ``A survey of evolutionary continuous dynamic optimization over two
  decades -- part {B},'' \emph{IEEE Transactions on Evolutionary Computation},
  vol.~25, no.~4, pp. 630--650, 2021.

\bibitem{kusakci2012constrained}
A.~O. Kusakci and M.~Can, ``Constrained optimization with evolutionary
  algorithms: a comprehensive review,'' \emph{Southeast Europe journal of soft
  computing}, vol.~1, no.~2, 2012.

\bibitem{omidvar2021reviewA}
M.~N. Omidvar, X.~Li, and X.~Yao, ``A review of population-based metaheuristics
  for large-scale black-box global optimization—{Part I},'' \emph{IEEE
  Transactions on Evolutionary Computation}, vol.~26, no.~5, pp. 802--822,
  2021.

\bibitem{omidvar2021reviewB}
------, ``A review of population-based metaheuristics for large-scale black-box
  global optimization—{Part II},'' \emph{IEEE Transactions on Evolutionary
  Computation}, vol.~26, no.~5, pp. 823--843, 2021.

\bibitem{li2016seeking}
X.~Li, M.~G. Epitropakis, K.~Deb, and A.~Engelbrecht, ``Seeking multiple
  solutions: An updated survey on niching methods and their applications,''
  \emph{IEEE Transactions on Evolutionary Computation}, vol.~21, no.~4, pp.
  518--538, 2016.

\bibitem{saini2021multi}
N.~Saini and S.~Saha, ``Multi-objective optimization techniques: a survey of
  the state-of-the-art and applications: Multi-objective optimization
  techniques,'' \emph{The European Physical Journal Special Topics}, vol. 230,
  no.~10, pp. 2319--2335, 2021.

\bibitem{yazdani2023GNBGgenerator}
\BIBentryALTinterwordspacing
D.~Yazdani, \emph{Generalized Numerical Benchmark Generator (GNBG)-Instance
  Generator (MATLAB Source Code)}, 2023 (accessed December 10, 2023). [Online].
  Available: \url{https://github.com/Danial-Yazdani/GNBG-Generator}
\BIBentrySTDinterwordspacing

\bibitem{hansen2009real}
N.~Hansen, S.~Finck, R.~Ros, and A.~Auger, ``Real-parameter black-box
  optimization benchmarking 2009: Noiseless functions definitions,'' Tech.
  Rep., 2009.

\bibitem{yazdani2018thesis}
D.~Yazdani, ``Particle swarm optimization for dynamically changing environments
  with particular focus on scalability and switching cost,'' Ph.D.
  dissertation, Liverpool John Moores University, Liverpool, UK, 2018.

\bibitem{yazdani2019scaling}
D.~Yazdani, M.~N. Omidvar, J.~Branke, T.~T. Nguyen, and X.~Yao, ``Scaling up
  dynamic optimization problems: A divide-and-conquer approach,'' \emph{IEEE
  Transactions on Evolutionary Computation}, vol.~24, no.~1, pp. 1--15, 2020.

\end{thebibliography}
\bibliographystyle{IEEEtran}

\end{document}